\documentclass[12pt]{amsart}
\newtheorem{theorem}{Theorem}[section]
\newtheorem{lemma}[theorem]{Lemma}
\newtheorem{corollary}[theorem]{Corollary}
\newtheorem{proposition}[theorem]{Proposition}
\theoremstyle{definition}
\newtheorem{definition}[theorem]{Definition}

\theoremstyle{remark}

\newcommand{\gap}{\hspace*{0.1em}}
\newcommand{\inv}{^{-1}}			

\newcommand{\sdpr}{\times\hspace{-0.088cm}\setlength{\unitlength}{1cm}\begin{picture}(0.5,0.0)\put(-0.075,0.01){\line(0,1){0.186}}\end{picture}\hspace{-0.5cm}}
\newcommand{\N}{\Bbb{N}}
\newcommand{\Z}{\Bbb{Z}}
\newcommand{\Q}{\Bbb{Q}}
\newcommand{\R}{\Bbb{R}}
\newcommand{\Aut}{{\rm Aut}}
\newcommand{\End}{{\rm End}}

\newcommand{\supp}{{\rm supp \hspace*{0.1em}}}
\newcommand{\Co}{{\mathcal C}}
\newcommand{\G}{{\mathcal G}}
\newcommand{\Null}{{\mathcal O}}
\newcommand{\U}{{\mathcal U}}
\renewcommand{\R}{{\mathcal R}}

\newcommand{\od}{{\overline d}}
\newcommand{\og}{{\overline g}}
\newcommand{\oh}{{\overline h}}

\newcommand{\om}{{\overline m}}

\newcommand{\ox}{{\overline x}}
\newcommand{\oy}{{\overline y}}
\newcommand{\oC}{{\overline C}}

\newcommand{\fC}{{\frak C}}
\newcommand{\ofC}{{\frak{\overline C}}}
\newcommand{\oalpha}{{\overline \alpha}}
\newcommand{\oeta}{{\overline \eta}}
\newcommand{\ta}{{\tilde a}}

\newcommand{\td}{{\tilde d}}
\newcommand{\tg}{{\tilde g}}
\newcommand{\tth}{{\tilde h}}

\newcommand{\tC}{{\tilde C}}
\newcommand{\tfC}{{\tilde{\frak C}}}

\begin{document}
\noindent
%
%
\title{\rm Free Division Rings of Fractions of Crossed Products}
\author{\vspace{0.9cm} Joachim Gr\"ater\medskip\\ 
}
\address{Universit\"at Potsdam, Institut f\"ur Mathematik, Karl-Liebknecht-Stra\ss e 24-25, 
D-14476 Potsdam OT Golm\\ Germany}
\email{graeter@uni-potsdam.de}
\begin{abstract}
Let $D$ be a division ring of fractions of a crossed product $F[G,\eta,\alpha]$ where $F$ is a skew field 
and $G$ is a group with Conradian left-order $\leq$. For $D$ we introduce the notion of freeness with 
respect to $\leq$ and show that $D$ is free in this sense if and only if $D$ can canonically be embedded 
into the endomorphism ring of the right $F$-vector space $F((G))$ of all formal power series in $G$ over 
$F$ with respect to $\leq$. From this we obtain that all division rings of fractions of $F[G,\eta,\alpha]$ 
which are free with respect to at least one Conradian left-order of $G$ are isomorphic and that they are 
free with respect to any Conradian left-order of $G$. Moreover, $F[G,\eta,\alpha]$ possesses a division 
ring of fraction which is free in this sense if and only if the rational closure of $F[G,\eta,\alpha]$ in 
the endomorphism ring of the corresponding right $F$-vector space $F((G))$ is a skew field.
\end{abstract}
\renewcommand{\thefootnote}{}
\footnote{
\smallskip\\ \footnotesize Date:\ \it \today}
\maketitle
\vspace*{-9.0cm}\ \\
\centerline{\rm OF GROUPS WITH CONRADIAN LEFT-ORDERS}
\vspace*{6.6cm}\ \\
\section{Introduction}\label{Introduction}

This paper deals with crossed products $F[G,\eta,\alpha]$ of a group $G$ with Conradian left-order 
$\leq$ over a skew field $F$. We are mainly interested in the case when $F[G,\eta,\alpha]$ possesses 
a division ring of fractions $D$. By this we mean a skew field $D$ which contains $F[G,\eta,\alpha]$ 
as a subring such that $D$ is generated by $F[G,\eta,\alpha]$ as a division ring. In Section 
\ref{Freie Quotientenschiefkoerper} we introduce the notion of freeness with respect to a Conradian 
left-order $\leq$ of $G$. Our first main result states that if $D$ meets this freeness condition then 
any non-zero $x \in D$ which is not a unit in $F[G,\eta,\alpha]$ possesses a "unique" normal form 
representation which can roughly be described as a local formal power series expansion over an archimedean 
ordered group where the coefficients are stricly simpler than $x$ with respect to the complexity discussed 
in \cite{D1,D2,DGH}. The archimedean ordered group corresponds to a factor group $C/C'$ where $(C',C)$ 
is a convex jump of $G$ with respect to the Conradian left-order $\leq$ and the coefficients are elements 
of the rational closure of $F[C',\eta,\alpha]$ in $D$. From this we derive our next main result that 
any element of $D$ can be comprehended as an endomorphism of the right $F$-vector space $F((G))$ of all 
formal power series in $G$ over $F$. In order to achieve this we first consider $F[G,\eta,\alpha]$ as a 
subring of the endomorphism ring $\End(F((G)))$ where the elements of $F[G,\eta,\alpha]$ act on $F((G))$ 
via multiplication from the left. Then Theorem \ref{Fortsetzungvonphi} states that $D$ is isomorphic 
to the rational closure $F(G,\eta,\alpha)$ of $F[G,\eta,\alpha]$ in $\End(F((G)))$. The following results 
are direct conclusions of our main theorems where the assumptions are as above. The rational closure 
$F(G,\eta,\alpha)$ of $F[G,\eta,\alpha]$ is a division ring if and only if $F[G,\eta,\alpha]$ possesses 
at least one division ring of fractions which is free with respect to $\leq$. In this case $F(G,\eta,\alpha)$ 
itself is free with respect to $\leq$, all such division rings of fractions are isomorphic to 
$F(G,\eta,\alpha)$, and every automorphisms of $F[G,\eta,\alpha]$ can be extended to an automorphism 
of $F(G,\eta,\alpha)$. For the division ring $D$ there are also other definitions of freeness which are 
due to I. Hughes and P.A. Linnell since any group with Conradian left-order is locally indicable 
and vice versa. In Section \ref{Quotientenschiefkoerper und Potenzreihen} we will study the relation 
between all these concepts of freeness and finally Section \ref{Hughes Theorems} contains a discussion 
of Hughes' two main theorems from \cite{H1,H2} in the light of our results.\smallskip\\
The topic of this paper is closely related to what has been investigated in \cite{J} where, however, 
the results are restricted to the special class of Conradian left-orders with maximal rank. 
Our approach does not require any further assumption. The theorems we prove hold true for all Conradian 
left-orders and they include the results from \cite{J} in this extended and generalized form.\smallskip\\
Throughout our investigation we make intensive use of two concepts that have originally been introduced 
by N.I. Dubrovin. Firstly, there is the complexity $cp(x)$ of an element $x$ in $D$ we already mentioned 
above. It is an ordinal which describes how $x$ is built from elements of $F[G,\eta,\alpha]$ by addition, 
multiplication, and the invers operation. It enables us to prove statements about the elements of $D$ 
by transfinite induction. Section \ref{Komplexitaet} provides all necessary background information 
which is needed to apply this concept. The proof of Theorem \ref{Fortsetzungvonphi} also requires 
Dubrovin's remarkable work on $v$-compatible and continuous endomorphisms of the right $F$-vector space 
$F((G))$ as given in \cite{D1,D2}.\\

\section{Conradian Groups}\label{Conrad}

\vspace*{0.3cm}
In order to make this paper self-contained we present in this short section all of the basics of 
left-ordered groups which will be used in the following. There will be nothing really new for the 
specialists but it hopefully improves the readability. A group $G$ equipped with a linear order $\leq$ 
is called left-ordered (with respect to $\leq$) if $a \leq b$ implies $ca \leq cb$ for all $a,b,c \in G$ 
and we also say that $\leq$ is a left-order of $G$. Right-orders are defined similarly but in this paper 
they are only needed on rare occasions. If they occur it will always be clear how to adopt the corresponding 
definitions and results stated previously for left-orders. For any left-order $\leq$ of $G$ the set 
$P_\leq =\{a \in G | e \leq a \}$ is called the positive cone of G with respect to $\leq$, where $e$ 
denotes the unit-element of $G$. Since $a \leq b \Longleftrightarrow a\inv b \in P_\leq$ for all 
$a,b \in G$ the left-order $\leq$ is uniquely determined by $P_\leq$. It is easily seen that 
$P_\leq \cup P_\leq \inv = G$, $P_\leq \cap P_\leq \inv = \{e\}$, and $P_\leq P_\leq \subseteq P_\leq$ where 
$A\inv = \{a\inv \mid a \in A\}$ and $AB = \{ab \mid a \in A, b\in B\}$ for arbitrary $A,B \subseteq G$ 
as usual. On the other hand if $P$ is a subset of $G$ satisfying $P \cup P\inv = G$, $P \cap P\inv = \{e\}$, 
and $P P \subseteq P$ then $a \leq' b \Longleftrightarrow a\inv b \in P$ for all $a,b \in G$ defines 
a left-order of $G$ such that $P$ is the corresponding positive cone. If $U$ is a subgroup of $G$ then 
any left-order $\leq$ of $G$ induces a left-order of $U$ (also denoted by $\leq$) and we will always 
consider $U$ as a left-ordered group in this way unless otherwise stated. For groups $G$ and $G'$ 
with left-orders $\leq$ and $\leq'$ respectively a group homomorphism $\varphi: G \longrightarrow G'$ is 
said to be order-preserving or an $o$-homomorphism if $a \leq b \Longrightarrow \varphi(a) \leq' \varphi(b)$ 
for all $a,b \in G$. We call $G$ and $G'$ order-isomorphic with respect to $\leq$ and $\leq'$ if there 
exists a group isomorphism $\varphi: G \longrightarrow G'$ which is order-preserving. The kernel of 
an $o$-homomorphism $\varphi: G \longrightarrow G'$ is a convex normal subgroup of $G$. On the other 
hand, if $N$ is a convex normal subgroup of $G$ then $P=\{aN | a \in P_\leq \}$ is a positive cone of 
the factor group $G/N$. The corresponding left-order $\leq'$ of $G/N$ is called the canonical left-order 
of $G/N$ with respect to $\leq$ and $\varphi: G \longrightarrow G/N$ is an $o$-homomorphism.

\begin{proposition}(cf. \cite{C})\label{convexohom}
Let $G$ and $G'$ be groups with left-orders $\leq$ and $\leq'$ respectively. If $\varphi: G \longrightarrow G'$ 
is a surjective $o$-homomorphism then the correspondence $C \longmapsto \varphi(C)$ induces an inclusion-preserving 
bijection between the set of all convex subgroups of $G$ containing the kernel of $\varphi$ and the set of 
all convex subgroups of $G'$.
\end{proposition}

\begin{definition}
Let $G$ be a group and let $\leq$ be a left-order of $G$. Then $G$ is called archimedean left-ordered 
and $\leq$ an archimedean left-order of $G$ if for all $a,b \in G$, $e < b$ there is an $n \in \N$ such 
that $a < b^n$.
\end{definition}

The following result is a generalization of the classical H\"older theorem (cf. \cite{Ho,C}):

\begin{theorem}\label{hoelderconrad}
Every archimedean left-ordered group is order-isomorphic to a subgroup of the naturally ordered 
additive group $\Bbb R$. Especially, any archimedean left-ordered group is commutative and possesses 
no proper convex subgroup.
\end{theorem}

Let $\leq$ be a left-order of the group $G$ and let $\Co_\leq$ denote the set of all convex subgroups 
of $G$. Then $\Co_\leq$ is a linearly ordered complete set with respect to the inclusion containing  
$\{e\}$ and $G$. A pair $(C',C)$ is called a convex jump in $G$ (with respect to $\leq$) whenever 
$C', C \in \Co_\leq$ and $C'$ is a lower neighbour of $C$, that is, $C'$ is a proper subgroup of $C$ 
and there is no convex subgroup of $G$ stricly between $C'$ and $C$ (cf. \cite{C,KM}).

\begin{definition}\label{cpluscminus}
Let $G$ be a left-ordered group with respect to $\leq$ and let $g$ be an arbitrary element of $G$. 
Then $C_g^+$ denotes the minimal convex subgroup of $G$ containing $g$. If $g$ is non-trivial then 
$C_g^-$ denotes the maximal convex subgroup of $G$ which does not contain $g$. 
\end{definition}

Cleary, $(C_g^-, C_g^+)$ is a convex jump in $G$ for all non-trivial $g$ in $G$. On the other hand if 
$(C',C)$ is a convex jump in $G$ then $C'=C_g^-$ and $C=C_g^+$ for any $g \in C$ which does not lie 
in $C'$. If in addition $C'$ is a normal subgroup in $C$ then the factor group $C/C'$ is called a 
factor of $\Co_\leq$.

\begin{definition}(cf. \cite{C,KM}). Let $G$ be a left-ordered group with respect to $\leq$. Then $G$ 
and $\leq$ are said to be Conradian if for any convex jump $(C',C)$ in $G$ with respect to $\leq$  the 
convex subgroup $C'$ is normal in $C$ such that the factor $C/C'$ is archimedean ordered with respect 
to the canonical order.
\end{definition}

Conradian left-orders can also be characterized without using convex subgroups (cf. \cite{C,KM}). In 
this respect it is worth mentioning that unlike us \cite{C} and \cite{KM} deal with right-orders 
(cf. \cite[Section 2]{GS}).

\begin{proposition}\label{conradkennzeichnung}
Let $G$ be a group with left-order $\leq$. Then the following statements are equivalent:
\begin{enumerate}
\item The left-order $\leq$ is Conradian.
\item For all $a,b \in G$ such that $e < a,b$ there is an $n \in \N$ satisfying $ab < (ba)^n$.
\item For all $a,b \in G$ such that $e < a < b$ there is an $n \in \N$ satisfying $b < a\inv b^n a$.
\item For all $a,b \in G$ such that $e < a,b$ there is an $n \in \N$ satisfying $a < b a^n$.
\end{enumerate}
\end{proposition}

As a consequence of this last result we obtain the first part of the following

\begin{proposition}\label{conraduntergruppe}
Let $G$ be a group with Conradian left-order $\leq$ and let $U$ be a subgroup of $G$. Then the left-order 
of $U$ induced by $\leq$ is Conradian. Moreover, if $(C',C)$ is a convex jump in $U$ then there 
exists a convex jump $(\oC', \oC)$ in $G$ extending $(C',C)$, that is, $\oC' \cap U = C'$ and $\oC \cap U = C$. 
Especially, $C/C'$ is isomorphic to a subgroup of $\oC/\oC'$.
\end{proposition}

{\bf Proof.} Clearly, $U$ is Conradian because of Proposition \ref{conradkennzeichnung}. Thus, let 
$(C',C)$ be a convex jump in $U$, let $\oC'$ be the maximal convex subgroup of $G$ such that 
$\oC' \cap U \subseteq C'$ and let $\oC$ be the minimal convex subgroup of $G$ containing $C$. Obviously, 
$\oC' \cap U$ and $\oC \cap U$ are convex subgroups of $U$ and $(\oC', \oC)$ is a convex jump in $G$. 
The canonical left-order of $\oC/\oC'$ induces on $(\oC \cap U)/(\oC' \cap U)$ the canonical left-order. 
Therefore, $C/(\oC' \cap U)$ and $C'/(\oC' \cap U)$ are two different convex subgroups of the archimedean 
ordered group $(\oC \cap U)/(\oC' \cap U)$ by Proposition \ref{convexohom}. Because of 
Theorem \ref{hoelderconrad} we conclude $\oC' \cap U = C'$ and $\oC \cap U = C$.\\
\qed\\

If $\leq$ is a Conradian left-order then Theorem \ref{hoelderconrad} shows that any factor of $\Co_\leq$ 
is order-isomorphic to a subgroup of the additive group $\Bbb R$.

\begin{definition}\label{defmaximalerrang} A Conradian left-order of a group $G$ has maximal rank if 
any factor of $\Co_\leq$ is order-isomorphic to a subgroup of the naturally ordered additive group $\Q$.
\end{definition}

This definition of maximality is equivalent to the one introduced in \cite{J}.

\begin{proposition}\label{abelschwohlgeordnet}
Any torsion-free abelian group $G$ possesses an order $\leq$ with maximal rank such that $\Co_\leq$ 
is well-ordered with respect to the inclusion.
\end{proposition}

{\bf Proof.} To simplify notation, let $G$ be written additively. We first assume that $G$ is divisible 
and regard $G$ as a $\Q$-vector space with basis $B=\{g_i \mid i \in I\}$ where $I$ is endowed with a 
well-order $\leq_I$. Any non-zero $a \in G$ can uniquely be written as 
$a = q_{i_1}g_{i_1} + \dots + q_{i_n}g_{i_n}$ with $n \in \N$, $i_1 <_I i_2 <_I \dots <_I i_n$, and 
non-zero rational numbers $q_{i_1}, \dots , q_{i_n}$. Defining $i_a := i_n \in I$ and $q(a):=q_{i_n}\in\Q$
it is easily seen that $P = \{ a \in G \mid a=0 \mbox{\ or\ } 0 < q(a)\}$ is a positive cone in $G$ 
with respect to an order $\leq$ of $G$ such that for any $i \in I$ 
\[H^-_i = \{a \in G \mid a = 0 \mbox{\ or\ } i_a < i\},\  
  H^+_i = \{a \in G \mid a = 0 \mbox{\ or\ } i_a \leq i\}\] 
are convex subgroups of $G$ such that $H^+_i/H^-_i$ and $\Q$ are isomorphic as ordered groups. Moreover, 
$H^+_i \subset H^+_j$ for all $i,j \in I, i <_I j$ and $H^+_{i(a)} = C^+_a$, $H^-_{i(a)} = C^-_a$ for 
all non-zero $a \in G$. This completes the proof if $G$ is divisible and we turn to the general situation. 
We consider the divisible hull $G'$ of $G$ for which the proposition has been proved and apply 
Proposition \ref{conraduntergruppe}.\\
\qed\\

\begin{proposition}\label{existenzmaximalerrang}
Let $G$ be a group with Conradian left-order $\leq$. Then there exists a Conradian left-order $\leq'$ 
with maximal rank such that $\Co_\leq \subseteq \Co_{\leq'}$. If $\Co_\leq$ is well-ordered with respect 
to the inclusion then $\leq'$ can be choosen that $\Co_{\leq'}$ is also well-ordered.
\end{proposition}

{\bf Proof.} If $(C',C)$ is a convex jump in $G$ with respect to $\leq$ then Proposition 
\ref{abelschwohlgeordnet} shows that there exists a left-order $\leq_{(C',C)}$ of $C/C'$ with maximal 
rank such that the corresponding convex subgroups are well-ordered with respect to the inclusion. 
Now, $\Co_\leq$ is a linearly ordered, complete, and subnormal system of subgroups of $G$ with respect 
to the inclusion and we can apply \cite[Proposition 3.2.1]{KM}. Thus, there exists a left-order $\leq'$ 
of $G$ such that $\Co_\leq \subseteq \Co_{\leq'}$ and for any convex jump $(C',C)$ in $G$ with respect 
to $\leq$ the canonical left-order of $C/C'$ with respect to $\leq'$ coincides with $\leq_{(C',C)}$. If 
$(\oC',\oC)$ is a convex jump in $G$ with respect to $\leq'$ then there is a convex jump $(C',C)$ in 
$G$ with respect to $\leq$ satisfying $C' \subseteq \oC' \subset \oC \subseteq C$ and $\oC'$ is normal 
in $\oC$. We conclude by Proposition \ref{convexohom} that $(\oC'/C', \oC/C')$ is a convex jump in $C/C'$ 
with respect to $\leq_{(C',C)}$ and therefore $\oC/\oC'$ is order-isomorphic to a subgroup of the additive 
group $\Q$. Finally, if $\Co_\leq$ well-ordered and if $C_1 \supset C_2 \supset \dots$ is a strictly 
decreasing chain of convex subgroups of $G$ with respect to $\leq'$ then there exists a convex jump 
$(C',C)$ of $G$ with respect to $\leq$ such that $C \supseteq C_k \supset C_{k+1} \supset \dots \supset C'$ 
for some $k \in \N$ since $\Co_\leq$ is well-ordered. But then 
$C/C' \supseteq C_k/C' \supset C_{k+1}/C' \supset \dots \supset C'/C'$ is a strictly decreasing chain 
of convex subgroups of $C/C'$ with respect to $\leq_{(C',C)}$ which is a contradiction.\\
\qed\\

\cite{G} provides a more detailed account of the results above and their proofs. Conradian left-orders 
can also be characterized by purely group-theoretic properties, that is, without mentioning terms from 
order theory.

\begin{definition}
Let $G$ be a group. Then $G$ is said to be indicable if there exists a surjective group homomorphism 
$\varphi: G \longrightarrow \Z$. And $G$ is called locally indicable whenever any non-trivial finitely 
generated subgroup of $G$ is indicable.
\end{definition}

According to the definition above, a group $G$ is indicable if and only if $G$ possesses a normal subgroup 
$N$ such that the factor group $G/N$ is an infinite cyclic group. Locally indicable groups have been 
introduced by Graham Higman in \cite{Hi} where he investigated zero-divisors and units in certain group 
rings. They also form an important class of groups in connection with the solvability of equations over 
groups (cf. \cite{Br}) and they occur naturally and in different ways in topology and geometry as 
fundamental groups of surfaces and manifolds, for instance as knot groups of classical knots. 
(cf. \cite{CR}).

\begin{theorem}\label{conradistlokalindizierbar}
Any Conradian group is locally indicable.
\end{theorem}

A more general version of this theorem has been proved by R.G. Burns and V.W. Hale and later on also 
by other authors (cf. \cite{BH}). A nice proof of Theorem \ref{conradistlokalindizierbar} is given in 
\cite[Section 10.1]{CR} which can be sketched as follows. Let $G=\langle g_1,\dots,g_n\rangle$ be a 
non-trivial finitely generated subgroup of a Conradian group. By Proposition \ref{conraduntergruppe} 
the group $G$ is also Conradian  and there exists $g \in \{g_1,\dots,g_n\}$ such that $G=C^+_g$. Then 
$N:=C^-_g$ is normal in $G$ and $G/N$ is a non-trivial, finitely generated, and torsion free abelian 
group. Thus, $G \cong \oplus_{i=1}^k \Z$ for some $k \in \N$ and now it is clear how to obtain a 
surjective group homomorphism $\varphi: G \longrightarrow G/N \longrightarrow \Z$.\smallskip\\

Even the convers of Theorem \ref{conradistlokalindizierbar} holds true:

\begin{theorem}\label{lokalindizierbaristconrad}
Any locally indicable group is Conradian.
\end{theorem}

Unlike Theorem \ref{conradistlokalindizierbar} a proof of Theorem \ref{lokalindizierbaristconrad} 
requires much more effort (cf. \cite{Br,RR,N}). The fact that the class of all Conradian groups coincides 
with the class of all locally indicable groups is of particular interest for this paper.

\section{Crossed Products}\label{Verschraenkte Produkte}

\vspace*{0.3cm}

In what follows, all rings are associative but not necessarily commutative, with unit element $1 \not= 0$ 
which is preserved by homomorphisms and inherited by subrings. If $R$ is a ring then we write $R^\times$ 
for the multiplicative group of units of $R$ and $R$ is called an integral domain if $R$ contains no 
non-trivial zero-divisor. By a skew field or division ring we shall understand a ring in which any 
non-zero element is a unit. In this section we summarize all definitions, propositions, and theorems 
about crossed products which are needed in this paper. We mainly refer to Section 1.1 of D.S. Passman's 
book {\it Infinite Crossed Products} (cf. \cite{Pa}).\smallskip\\
Let $R$ be a ring and $G$ a group. A crossed product of $G$ over $R$ is a ring $S$ containing $R$ as 
a subring which turns $S$ naturally into a left $R$-module which is free with basis 
$X_G = \{x_g \ |\ g \in G\}$ where $x_g\not=x_h$ if $g\not=h$ for all $g,h \in G$. By this we mean that 
the additive structure of $S$ as a ring coincides with the additive structure of $S$ as an $R$-module and 
with respect to the multiplication there are two mappings $\eta: G \times G \longrightarrow R^\times$ 
and $\alpha: G \longrightarrow \Aut R, g \longmapsto \alpha_g$ having the following properties:
\begin{enumerate}
\item [(1)]$1x_e$ is the unit element of $S$,
\item [(2)]$ax_e \cdot bx_e = abx_e$ for all $a,b \in R$,
\item [(3)]$ax_e\cdot 1x_g = ax_g$ for all $a\in R$ and $g\in G$,
\item [(4)]$1x_g \cdot 1x_h = \eta(g,h)x_{gh}$ for all $g,h\in G$,
\item [(5)]$1x_g \cdot ax_e = \alpha_g(a)x_g$ for all $a\in R$ and $g\in G$.
\end{enumerate}

From this we obtain the following equations for all $g,h,l \in G$ and $a \in R$:

\begin{enumerate}
\item[(6)] $\eta(g,e)=\eta(e,g) = 1$,
\item[(7)] $\alpha_g\alpha_h(a) = \eta(g,h)\alpha_{gh}(a) \eta(g,h)^{-1}$,
\item[(8)] $\alpha_g(\eta(h,l))\eta(g,hl) = \eta(g,h)\eta(gh,l)$.
\end{enumerate}

Because of $(1)$ and $(2)$ we usually write $a$ instead of $ax_e$ for all $a \in R$ and $x_g$ instead 
of $1x_g$ for all $g \in G$ such that $(3)$ provides $a \cdot x_g = ax_e \cdot 1x_g = ax_g$. 
Then any element $x \in S$ can uniquely be written as a finite sum $x = \sum a_gx_g$ with $a_g \in R$, 
$g \in G$ and the support of $x$ is defined as the set of all $g \in G$ with non-zero $a_g$. Moreover, 
the multiplication in $S$ turns into the following form:
\[(\sum a_g x_g) \cdot (\sum b_h x_h) = \sum a_g x_g \cdot b_h x_h \ \mbox{where}\ 
a_g x_g \cdot b_h x_h = a_g\alpha_g(b_h)\eta(g,h)x_{gh}.\]
This shows that the structure of $S$ as a ring and as an $R$-module is uniquely determined by $R,G,\eta$, 
and $\alpha$ and therefore we also write $R[G,\eta,\alpha]$ instead of $S$.\smallskip\\
On the other hand if a ring $R$ and a group $G$ are given together with two mappings   
$\eta: G \times G \longrightarrow R^\times$ and 
$\alpha: G \longrightarrow \Aut R, g \longmapsto \alpha_g$ such that (6),(7), and (8) are satisfied 
then there exists a crossed product $S=R[G,\eta,\alpha]$ of $G$ over $R$.\smallskip\\
Any $x_g, g\in G$ is a unit in $R[G,\eta,\alpha]$ such that $x_g a x\inv_g = \alpha_g(a)$ for all $a \in R$. 
This shows that the conjugation by $x_g$ defines a ring automorphism of $R[G,\eta,\alpha]$ which extends 
$\alpha_g$ to $R[G,\eta,\alpha]$ and which will be denoted by $\alpha_g$ again, that is,
\[\alpha_g(\sum a_h x_h) = \sum \alpha_g(a_h) x_gx_hx_g\inv.\]
If $T$ is a subring of $R$ and $U$ a subgroup of $G$ such that $\alpha_g(T)=T$ for all $g \in U$ and 
and $\eta(g,h) \in T^\times$ for all $g,h \in U$ then the subring of $R[G,\eta,\alpha]$ which is generated 
by $T$ and $X_U$ is a crossed product of $U$ over $T$ in an obvious way. By abuse of 
notation, this ring will also be denoted by $T[U,\eta,\alpha]$. The proof of the following proposition 
is straightforward.

\begin{proposition}\label{RGverschraenktueberRN}
With the notation as above the following holds true: If $N$ is a normal subgroup of $G$ and if $\frak G$ 
is a transversal for $N$ in $G$ such that $e \in \frak G$ then $R[G,\eta,\alpha]$ is a crossed product 
of $G/N$ over $R[N,\eta,\alpha]$ with basis $X_{G/N}=\{x_g\ |\ g\in\frak G\}$. Moreover, for any $g \in G$ 
the restriction $\oalpha_g$ of $\alpha_g$ to $R[N,\eta,\alpha]$ is a ring automorphism of 
$R[N,\eta,\alpha]$.
\end{proposition}

If $R[G,\eta,\alpha]$ is a crossed product of $G$ over $R$ and $S$ a ring containing $R$ as a subring 
such that any $\alpha_g, g\in G$ can be extended uniquely to a ring automorphism $\oalpha_g$ of $S$ 
then conditions $(6) - (8)$  apply accordingly for $\eta$ and 
$\oalpha: G \longrightarrow \Aut S, g \longmapsto \oalpha_g$. This is rather obvious for property (6) 
and (8). In order to prove property (7) we first observe that the automorphism $\alpha_h$ can be extended 
to $S$ by $\oalpha_h$ on the one hand, and by the automorphism of $S$ arising from the composition of 
$\oalpha_{gh}$, the conjugation by $\eta(g,h)$, and $\oalpha\inv_g$ on the other. Because of the 
uniqueness of the extension we are done. Thus we have proved

\begin{proposition}\label{verschraenktesproduktueberd}
With the notation as above the following holds true: Let $R$ be a subring of the ring $S$ and let 
$R[G,\eta,\alpha]$ be a crossed product of $G$ over $R$. If for any $g \in G$ there exists a unique 
ring automorphism $\oalpha_g$ of $S$ extending $\alpha_g$ then $R[G,\eta,\alpha]$ coincides in an 
obvious way with the subring of $S[G,\eta,\oalpha]$ which is generated by $R$ and $X_G$.
\end{proposition}

The following proposition has been used before by other authors in different variants. For the sake 
of completeness, we will briefly outline a proof for crossed products of groups which are left-ordered 
(cf. \cite{Bo}).

\begin{proposition}\label{Einheitengruppe}
Let $R=F[G,\eta,\alpha]$ be a crossed product of a group $G$ with left-order $\leq$ over a skew field 
$F$. Then $R$ is an integral domain with group of units $F^\times X_G$. If 
$\varphi: R \longrightarrow R$ is a ring automorphism then $\varphi(F) = F$, $\varphi(FX_G) = FX_G$, 
and $R$ is a crossed product of $G$ over $F$ with respect to the basis 
$X'_G = \{x'_g \ |\ g \in G\}$ where $x'_g = \varphi(x_g)$ for all $g \in G$.
\end{proposition}

{\bf Proof.} Let $x,y$ be non-zero elements in $R$ where $x = a_1x_{g_1} + \dots + a_mx_{g_m}$ and 
$y = b_1x_{h_1} + \dots + b_nx_{h_n}$ with $a_1,\dots,a_m,b_1,\dots,b_n \in F^\times$. We assume 
$g_1 < g_2 < \dots < g_m$ and $h_1 < h_2 < \dots < h_n$. Then $g_ih_1 < \dots < g_ih_n$ follows for 
all $i \in \{1,\dots,m\}$. Let $k$ be in $\{1,\dots,m\}$ such that $g_kh_n \geq g_1h_n, \dots,g_mh_n$. 
We conclude
\[g_kh_n > g_ih_j \mbox{\ for \ } (i,j) \not= (k,n).\]
Therefore, $a_k\alpha_{g_k}(b_n)\eta(g_k,h_n)$ is the non-zero coefficient of $g_kh_n$ in the 
presentation of $xy$ which proves the first claim of the proposition. By similar arguments 
there exists $l \in \{1,\dots,m\}$ such that
\[g_lh_1 < g_ih_j \mbox{\ for \ } (i,j) \not= (l,1)\]
and $g_lh_1$ is also an element of the support of $xy$. Now, let $x$ be a unit in $R$ and let $xy=1$. 
In order to show $x \in F^\times X_G$ we assume $1<m$. Then there exists an $i \in \{1,\dots,m\}$ 
different from $l$ and $g_lh_1 < g_ih_1 \leq g_kh_n$ follows. This shows that the support of 
$xy$ contains more than one element which contradicts $xy = 1$. Since any element of $F^\times X_G$ 
is obviously a unit in $R$ the second claim of the proposition is shown.\smallskip\\
Now, let $\varphi: R \longrightarrow R$ be a ring automorphism. Since $F^\times X_G$ is the group of 
units in $R$ we obtain $\varphi(F^\times X_G) = F^\times X_G$, that is $\varphi(FX_G) = FX_G$. 
Next we prove $\varphi(F) \subseteq F$. Then $\varphi^{-1}(F) \subseteq F$ follows and we get 
$\varphi(F) = F$ as claimed. Let $a$ be in $F$ different from $0$ and $1$. Since $a-1$ is a unit in 
$R$ there exist $b \in F^\times$ and $g \in G$ such that $\varphi(a) =bx_g$. Hence, 
$\varphi(a-1) = \varphi(a) - \varphi(1) = bx_g - x_e \in F^\times X_G$ which shows $g = e$. The claim 
that $R$ is a crossed product of $G$ over $F$ with respect to the basis $X'_G$ can now be derived by 
straightforward computation.\\
\qed\\

Before presenting the last result of this section we discuss two examples where $R$ is assumed to be a 
(left and right) Ore domain (cf. \cite[Section 5]{Co}). Firstly, let $S=R[\langle g \rangle, \eta, \alpha]$ 
be a crossed product of an infinite cyclic group $G= \langle g \rangle$ over $R$. For any $n\in \Z$ 
there exists a unit $a_n \in R^\times$ such that $x_{g^n} = a_n x^n_g$. Thus $\{x^n_g \mid n \in \Z\}$ 
is an $R$-basis of the free $R$-module $S$, that is, $S$ can be seen as the skew Laurent polynomial 
ring in $x_g$ over $R$ which shows that $R[\langle g \rangle, \eta, \alpha]$  is an Ore domain. In the 
second example $G$ is isomorphic to a subgroup of the additive group $\Q$. For given $a,b \in S$ there 
exist $g,g_1,\dots,g_n \in G$ such that $a,b \in Rx_{g_1} + \dots + Rx_{g_n}$ and 
$\langle g_1,\dots,g_n \rangle = \langle g \rangle$. Therefore, $a,b$ belong to the subring of $S$ 
generated by $R$, $x_g$, and $x_{g\inv}$, that is, $a,b \in R[\langle g \rangle, \eta, \alpha]$ which 
is an Ore domain as mentioned above. Thus, $R[G, \eta, \alpha]$ is an Ore domain too.

\begin{proposition}\label{conradore}
Let $R=F[G,\eta,\alpha]$ be a crossed product of a group $G$ over a skew field $F$. If $G$ possesses 
a Conradian left-order $\leq$ such that $\Co_\leq$ is well-ordered with respect to the inclusion then 
$R$ is an Ore domain.
\end{proposition}

{\bf Proof.} Because of Proposition \ref{existenzmaximalerrang} let $\leq$ be a Conradian left-order with 
maximal rank. We assume that $F[G,\eta,\alpha]$ is not an Ore domain and choose $C \in \Co_\leq$ 
minimal with respect to the inclusion such that $F[C,\eta,\alpha]$ is not an Ore domain. Then $C$ must 
have a lower neighbour $C' \in \Co_\leq$, that is, $(C',C)$ is a convex jump. According to Proposition 
\ref{RGverschraenktueberRN}, $F[C,\eta,\alpha]$ is a crossed product of $C/C'$ over $F[C',\eta,\alpha]$. 
Because of the minimality of $C$ we conclude that $F[C',\eta,\alpha]$ is an Ore domain. Since $C/C'$ is 
isomorphic to a subgroup of the additive group $\Q$ the second example above provides the desired 
contradiction.\\
\qed\\

\section{The Complexity}\label{Komplexitaet}

\vspace*{0.3cm}

If $D$ is a division ring with subring $R$ then the rational closure of $R$ in $D$ is the minimal 
subdivision ring of $D$ containing $R$. If $D$ itself is this rational closure then $D$ is called 
a division ring of fractions of $R$. In this case any element $d$ of $D$ is built from 
elements of $R$ by means of addition, multiplication, and the invers operation, that is, $d$ possesses 
a presentation as a rational expression of elements from $R$. Of course, such a presentation is not 
unique and in general it does not even exist a unique normal form in which any presentation can be 
transformed. Despite this difficulty it is shown in \cite{DGH} that it is always possible to assign 
ordinals to the elements of a division ring of fractions in order to describe the complexity 
of the expressions which are needed to present an element of $D$ by elements from $R$. 
For $d \in D$ this ordinal is called the complexity of $d$ denoted by $cp(d)$. The complexity allows 
us to prove statements about the elements of $D$ by transfinite induction which is an essential tool  
for proving our main theorems. In this section we give a short introduction into the basic concepts of 
the complexity and we explain how to apply the results from \cite{DGH}. We also present some results 
which are not included in \cite{DGH} but which are needed here.\smallskip\\
In this section $T$ always denotes a ring (which is denoted by $R$ in \cite{DGH}). If $A$ is a subset 
of $T$ then $A\inv$ consists of all $a\inv$ where $a \in A \cap T^\times$. By the rational closure 
of $A$ in $T$ we will understand the smallest subring $D$ of $T$ which contains $A$ such that 
$D\inv \subseteq D$. We fix a subset $\G$ of $T^\times$ (which is denoted by $S$ in \cite{DGH}) and 
assume that $\G$ is a subgroup of $T^\times$ and $-g \in \G$ for all $g \in \G$, that is $-\G = \G$. 
Finally, $D$ denotes the rational closure of $\G$ in $T$. In our applications $T$ always contains a 
subring $R = F[G,\eta,\alpha]$ which is a crossed product of a left-ordered group $G$ over a skew field 
$F$. Then $\G$ will be the group of units of $R$, that is, $\G = F^\times X_G$ and $R$ has no non-trivial 
zero-divisor. If $D$ is a division ring then $D$ is a division ring of fractions of $F[G,\eta,\alpha]$.
\smallskip\\
As mentioned above any $a \in D$ is assigned an ordinal denoted by $cp(a)$, e.g., $cp(0)=0$ and 
$cp(a) = 1$ if and only if $a \in \G$. In case of a crossed product as described above we obtain 
$cp(a)=n$ for any $a = a_1x_{g_1} + \dots + a_nx_{g_n} \in R$ if $a_1,\dots,a_n \in F^\times$ and if 
$g_1,\dots,g_n \in G$ are pairwise distinct. For further applications it is not important to know precisely 
how $cp(a)$ is defined. It will be enough to know some basic rules which will enable us to apply results 
from \cite{DGH} or to derive further properties. For $a,b \in D$ we write $a \unlhd b$  if $cp(a) \leq cp(b)$ 
and $a \lhd b$ if $cp(a) < cp(b)$. In the first case we say that $a$ is simpler than $b$ and in 
the second that $a$ is strictly simpler than $b$.\smallskip\\
Besides the complexity of single elements we also need their formal sums. In order to explain this we 
first have to fix an ordinal $\Delta$ satisfying $cp(a) \leq \Delta$ for all $a \in D$ and consider 
the free commutative monoid $\N(\Delta)$ freely generated by $\{ \alpha \ |\  0 < \alpha \leq \Delta\}$. 
In Section 2 of \cite{DGH} a well-ordering $\leq$ of $\N(\Delta)$ is defined which is compatible with 
the addition of $\N(\Delta)$. For different $x,y$ in $\N(\Delta)$ the relation $x < y$ can be explained 
as follows: There exist $\alpha_1,\dots,\alpha_k \leq \Delta$ with $\alpha_1 > \dots > \alpha_k$ and 
$n_1,\dots,n_k,m_1,\dots,m_k \in \N \cup \{0\}$ such that $x = n_1\alpha_1 + \dots + n_k\alpha_k$ and 
$y = m_1\alpha_1 + \dots + m_k\alpha_k$. Since $x \not= y$ there is a minimal $i \in \{1, \dots ,k\}$ 
satisfying $n_i \not= m_i$. Then $x < y$ if and only if $n_i < m_i$.\smallskip\\
{\it Additive Decompositions.} If an element $a$ of $D$ can be written as a sum of elements of $D$ which 
are all strictly simpler than $a$ then $a$ is said to be additively decomposable (otherwise $a$ is called 
additively indecomposable). Among all these decompositions there are some which are minimal in a certain 
sense. These specific decompositions can be described by a function $\tau : D \longrightarrow \N(\Delta)$ 
introduced in \cite{DGH} as follows. Let $a$ be additively decomposable and let $a_1,\dots,a_n$ be non-zero 
elements of $D$ such that $a =a_1 + \dots + a_n$ and $a_1,\dots,a_n \lhd a$. Then $a =a_1 + \dots + a_n$ is 
called a complete additive decomposition of $a$ if $cp(a_1) + \dots + cp(a_n) \in \N(\Delta)$ is minimal 
in a certain sense with respect to the well-order $\leq$ of $\N(\Delta)$ introduced above. Following 
\cite{DGH} we define $\tau(a) := cp(a_1) + \dots + cp(a_n)$ whenever $a =a_1 + \dots + a_n$ is a complete 
additive decomposition. For formal reasons it is useful to define $\tau(a)$ even if $a$ is additively 
indecomposable. In this case we put $\tau(a) = cp(a)$ and call $a$ itself a complete additive decomposition 
of $a$. Thus, $a \in D$ is additively indecomposable if and only if $\tau(a) = cp(a)$.\smallskip\\
{\it Multiplicative Decompositions.} Similar to $\tau$ there exists a function $\rho$ which describes 
minimal multiplicative decompositions of elements in $D$. Let $a \in D$ be additively indecomposable. 
Then $a$ is said to be multiplicatively decomposable if $a$ can be written as a product of elements 
of $D$ which are all strictly simpler than $a$. Otherwise $a$ is called an atom. For example, $0$ and 
the elements from $\G$ are atoms. All other atoms are called proper. Formulations like {\it $a$ is 
multiplicatively decomposable} or {\it a is an atom} always include that $a$ is additively indecomposable. 
If $a$ is multiplicatively decomposable then a decomposition $a = a_1 \cdot \dots \cdot a_n$ of $a$ 
such that $a_1,\dots,a_n \lhd a$ is called complete if $cp(a_1) + \dots + cp(a_n) \in \N(\Delta)$ is 
minimal in a certain sense with respect to the well-order $\leq$ of $\N(\Delta)$ introduced above. 
In this case we define $\rho(a) = cp(a_1) + \dots + cp(a_n)$ and moreover if $a$ is an atom we put 
$\rho(a) = cp(a)$ and call $a$ itself a complete multiplicative decomposition of $a$. Thus, if $a$ 
is additively indecomposable then $a$ is an atom if and only if $\rho(a) = cp(a)$.\smallskip\\
We return to the example where $T$ contains the subring $R=F[G,\eta,\alpha]$ with a skew field $F$, a 
left-ordered group $G$, and $\G = F^\times X_G$. Clearly, the subring $D_{[0]}$ of $T$ generated by 
$\G$ coincides with $F[G,\eta,\alpha]$. Let $a = a_1x_{g_1} + \dots + a_nx_{g_n} \in R$ such 
that $a_1,\dots,a_n \in F^\times$ and $g_1,\dots,g_n \in G$ are pairwise distinct. Then, 
$a = a_1x_{g_1} + \dots + a_nx_{g_n}$ is a complete additive decomposition of $a$. Moreover, $a$ is 
no unit in $R$ if $n>1$ by Proposition \ref{Einheitengruppe}. If $n=2$ and if $a$ is a unit in $T$ 
then $a\inv = (a_1x_{g_1} + a_2x_{g_2})\inv$ is a proper atom in $D$.\smallskip\\ 
The next two propositions will be needed in Section \ref{Quotientenschiefkoerper und Potenzreihen}. 
They are not included in \cite{DGH} but they can easily be derived.\smallskip

\begin{proposition} \label{neuezerlegungen}
Let $x \in D$, $cp(x) > 1$ and $g \in \G$. If $x$ is additively indecomposable then $gx$ is also 
additively indecomposable and if $x$ is an atom then $gx$ is also an atom. If $x$ is additively 
indecomposable and $x = x_1 \cdot x_2 \cdot\dots\cdot x_n$ is a complete multiplicative decomposition 
of $x$ then $cp(x_i) > 1$ for all $i=1,\dots,n$ and $gx = (gx_1) \cdot x_2 \cdot\dots\cdot x_n$ is a 
complete multiplicative decomposition of $gx$. If $x$ is additively decomposable with complete additive 
decomposition $x = x_1 + \dots + x_n$  then $gx = gx_1 + \dots + gx_n$ is a complete additive decomposition 
of $gx$.
\end{proposition}

{\bf Proof.} Let $x$ be additively indecomposable and assume that $gx$ is additively decomposable. 
Then $gx$ can be written as $gx = a_1 + a_2$ where $a_1, a_2 \lhd gx$ and therefore 
$x = g^{-1}a_1 + g^{-1}a_2$ where $cp(x) = cp(gx)$ and $cp(a_i) = cp(g^{-1}a_i)$ for  $i=1,2$ by 
\cite[Proposition 4.8]{DGH}. Thus, $x$ is additively decomposable - a contradiction. Next, let $x$ 
be an atom. Again, we assume that $gx$ is not an atom. Since $gx$ is additively indecomposable there 
exist $a_1,a_2 \in D$ such that $gx = a_1 a_2$ where $a_1, a_2 \lhd gx$. Then $x = (g^{-1}a_1) a_2$ 
with $g^{-1}a_1, a_2 \lhd x$ which is a contradiction. Now, let $x$ be additively indecomposable and 
let $x = x_1 \cdot x_2 \cdot\dots\cdot x_n$ be a complete multiplicative decomposition of $x$. Any 
$x_i$, $i=1,\dots,n$ is an atom by \cite[Theorem 4.6]{DGH}. Without loss of generality we just 
show $cp(x_1) > 1$ and assume $cp(x_1) = 1$. Then $n>1$ and $x'=x_1x_2$ is additively indecomposable 
because of \cite[Lemma 4.5]{DGH}. Furthermore, $x'=x_1x_2$ is a complete multiplicative decomposition 
of $x'$ by \cite[Theorem 4.6]{DGH} such that $\rho(x')=cp(x_1) + cp(x_2)$. On the other side $x'$ is 
an atom since $x_1 \in \G$ by assumption. This leads to the contradiction $\rho(x') = cp(x')$. 
Furthermore, let $gx = a_1 \cdot a_2 \cdot\dots\cdot a_m$ be a complete multiplicative decomposition 
of $gx$. According to \cite[Proposition 4.2,4.8]{DGH} and \cite[Theorem 4.4]{DGH} we obtain
\begin{eqnarray*}
\rho(x) &\leq & cp(g^{-1}a_1) + cp(a_2) + \dots + cp(a_m)\\
        & = &   cp(a_1) + cp(a_2) + \dots + cp(a_m) = \rho(gx)\\
        &\leq & cp(g x_1) + cp(x_2) + \dots + cp(x_n)\\
        & = &   cp(x_1) + cp(x_2) + \dots + cp(x_n) = \rho(x)
\end{eqnarray*}
where $\rho$ is the function introduced above. Thus, $\rho(gx) = cp(g x_1) + cp(x_2) + \dots + cp(x_n)$ 
which implies the claim because of \cite[Proposition 4.2]{DGH}.\smallskip\\
Finally, let $x$ be additively decomposable with complete additive decomposition $x = x_1 + \dots + x_n$ 
and let $gx = a_1 + \dots + a_m$ be a complete additive decomposition of $gx$. We apply 
\cite[Proposition 3.2,4.8]{DGH} and \cite[Theorem 3.4]{DGH} and obtain by arguments similar to those 
we just used above
\begin{eqnarray*}
\tau(x) &\leq & cp(g^{-1}a_1) + \dots + cp(g^{-1}a_m) = cp(a_1) + \dots + cp(a_m) = \tau(gx)\\
        &\leq & cp(g x_1) + \dots + cp(gx_n) = cp(x_1) + \dots + cp(x_n) = \tau(x).
\end{eqnarray*}
Therefore, $\tau(gx) = cp(gx_1) + \dots + cp(gx_n)$ which implies the claim because of 
\cite[Proposition 3.2]{DGH}.\\
\qed\\

{\bf Remark.} Corresponding statements hold true for terms like $xg$.

\begin{definition}\label{lambdausw}
Let the notation be as above.
\begin{enumerate}
\item $\Lambda = \{ cp(x) \mid x \in D, x \not=0, \mbox{\ and $x$ is additively indecomposable}\}$.
\item If $\lambda \in \Lambda$ then $D_{\leq \lambda}$ denotes the subgroup of the additive group of 
      $D$ which is generated by all $a \in D$ satisfying $cp(a) \leq \lambda$.
\item If $\lambda \in \Lambda$ then $D_{< \lambda}$ denotes the subgroup of the additive group of $D$ 
      which is generated by all $a \in D$ satisfying $cp(a) < \lambda$.
\end{enumerate}
\end{definition}

\begin{proposition}\label{Dleqlambda}
With the notation as above the following hold true:
\begin{enumerate}
\item For all $\lambda \in \Lambda$, $a \in D_{\leq\lambda}$ there are $a_1,\dots,a_n \in D$ such 
      that $a = a_1 + \dots + a_n$ and $cp(a_1),\dots,cp(a_n) \leq \lambda$.
\item $D_{<\lambda} = \{a \in D \ |\ cp(a) < \lambda\}$ for all $\lambda \in \Lambda$.
\item Let $a \in D, a \not=0$, and $\lambda \in \Lambda$. If $a = a_1 + \dots + a_n$ is a complete 
      additive decomposition of $a$ then $a$ belongs to $D_{\leq \lambda}$ if and only if 
      $cp(a_1),\dots,cp(a_n) \leq \lambda$.
\item If $\lambda \in \Lambda$, $a \in D_{\leq \lambda}$, and $b \in D$ such that $b \lhd a$ then 
      $b \in D_{\leq \lambda}$ .
\item $D_{\leq \lambda} \subset D_{\leq \lambda'}$ for all $\lambda, \lambda' \in \Lambda$ such that 
      $\lambda < \lambda'$.
\item Let $a \in D$, $cp(a) > 1$ be additively indecomposable. For all $a_1,\dots,a_n \in D$ such that 
      $a_1,\dots,a_n \lhd a$ there exists $\lambda \in \Lambda$ with $\lambda < cp(a)$ and 
      $a_1,\dots,a_n \in D_{\leq \lambda}$.
\item For all $\lambda \in \Lambda$ and non-zero $x \in D_{< \lambda}$ there exists $\lambda' \in \Lambda$ 
      such that $\lambda' < \lambda$ and $x \in D_{\leq \lambda'}$.
\end{enumerate}
\end{proposition}

{\bf Proof.} (1): Clearly, $-1 \in \G$ since $\G = -\G$ by assumption and \cite[Proposition 4.8]{DGH} 
proves the claim. 
\\
(2): This statement follows from \cite[Proposition 3.1]{DGH}.\\
(3): Let $a$ be non-zero and let $a = a_1+\dots+a_n$ be a complete additive decomposition of $a$. 
If $cp(a_1),\dots,cp(a_n) \leq\lambda$ then $a \in D_{\leq \lambda}$ by definition of 
$D_{\leq \lambda}$. Conversely, let $a$ be in $D_{\leq \lambda}$, that is, $a = c_1 + \dots + c_m$ 
where $c_1,\dots,c_m \in D$ and $cp(c_1),\dots,cp(c_m) \leq \lambda$. It will be enough to show that 
for all $i = 1,\dots,n$ there is an index $j = 1,\dots,m$ such that $a_i \unlhd c_j$. Let us assume 
that there exists an $i \in \{1,\dots,n\}$ such that $c_j \lhd a_i$ for all $j = 1,\dots,m$ which means 
$cp(c_1) + \dots + cp(c_m) < cp(a_1) + \dots + cp(a_n)$ with respect to the order defined in $\N(\Delta)$. 
Then $cp(c_1) + \dots + cp(c_m) < \tau(a)$ 
because of \cite[Proposition 3.2]{DGH}. Applying \cite[Lemma 3.3]{DGH} we conclude  
$a = c_1 + \dots + c_m \lhd a$, a contradiction.\\ 
(4): Let $a = a_1 + \dots + a_n$ and $b = b_1 + \dots + b_m$ be two complete additive decompositions 
of $a$ and $b$ respectively. Then $a \not= 0$ since $b \lhd a$ and (3) yields 
$cp(a_1), \dots ,cp(a_n) \leq \lambda$. We apply \cite[Theorem 3.4]{DGH} to obtain 
$\tau(b) \lhd \tau(a)$. Thus, for any $i = 1, \dots ,m$ there exists an index $j = 1, \dots ,n$ such 
that $b_i \unlhd a_j$, that is $cp(b_i) \leq cp(a_j) \leq \lambda$. This implies 
$b \in D_{\leq \lambda}$.\\
(5): Obviously, $\lambda < \lambda'$ provides $D_{\leq \lambda} \subseteq D_{\leq \lambda'}$. We choose 
a non-zero $a \in D$ which is additively indecomposable such that $cp(a) = \lambda'$. Then $a$ is a complete 
additive decomposition of $a$ and $a \in D_{\leq \lambda}$ cannot hold true because of (3).\\
(6): It will be enough to show that for any $a' \in D$ satisfying $a' \lhd a$ there exists 
$\lambda \in \Lambda$ such that $a' \in D_{\leq\lambda}$ and $\lambda < cp(a)$. If $a'=0$ we put 
$\lambda = 1$. Thus, let $a' \not = 0$. If $a'$ is additively indecomposable we define $\lambda = cp(a')$. 
Now, let $a'$ be additively decomposable and let $a' = a'_1 + \dots + a'_m$ be a complete additive 
decomposition of $a'$ where $m > 1$. We assume $cp(a) \leq cp(a'_i)$ for some $i \in \{1,\dots,m\}$. 
Then $cp(a) < cp(a'_1) + \dots + cp(a'_m)$ with respect to the given order of $\N(\Lambda)$ and therefore 
$cp(a'_1) + \dots + cp(a'_m) = \tau(a') < cp(a')$ by \cite[Proposition 3.2]{DGH}. This implies the 
contradiction $cp(a) < cp(a')$ and therefore $cp(a'_i) < cp(a)$ for all $i = 1,\dots,m$. We may choose 
$\lambda$ as the maximum of all $cp(a'_1),\dots, cp(a'_m)$ since any $a'_i$, $i=1,\dots,m$ is additively 
indecomposable by \cite[Theorem 3.6]{DGH}.\\
(7): If $x$ is additively indecomposable we choose $\lambda' = cp(x)$. Thus, let $x = a_1 + \dots + a_n$ 
be a complete additive decomposition of $x$ where $n > 1$ and let $a$ be a non-zero element of $D$ 
with $\lambda = cp(a)$. Because of  $a_1, \dots , a_n \lhd x \lhd a$ and (6) there exists 
$\lambda' \in \Lambda$ satisfying $\lambda' < \lambda$ and $a_1, \dots , a_n \in D_{\leq \lambda'}$, 
that is $x \in D_{\leq \lambda'}$.\\
\qed

\section{Free Division Rings of Fractions}\label{Freie Quotientenschiefkoerper}

\vspace*{0.4cm}

In this section we will use all definitions and results from Section \ref{Verschraenkte Produkte} 
without any further explanation and fix the following notation. $R=F[G,\eta,\alpha]$ always 
denotes a crossed  product of $G$ over $F$ with basis $X_G=\{x_g \mid g\in G\}$ where $G$ is a 
group with left-order $\leq$ and $F$ a skew field such that $\G = F^\times X_G$ is the group of 
units in $R$. The fact that $G$ is left-ordered is not always needed but this does not appear to be 
relevant here. Furthermore, $D$ is a division ring of fractions of $R$. If $U$ is a subgroup of $G$ 
then $F^\times X_U$ is a subgroup of $\G$ and $D_U$ denotes the rational closure of $F^\times X_U$ 
in $D$, that is, $D_U$ is a division ring of fractions of the crossed product $F[U,\eta,\alpha]$.
\smallskip\\
Let $U$ be a normal subgroup of $G$ and $h$ an element of $G$. If $t \in \G$ is a preimage of $h$ 
with respect to the canonical group homomorphism $\G \longrightarrow G, ax_g \longmapsto g$, that is, 
$t = ax_h$ for some $a\in F^\times$, then $tdt^{-1} \in D_U$ holds true for all $d \in D_U$. This 
can easily be shown for instance by transfinite induction on the complexity $cp(d)$ of $d$ since $U$ is 
normal in $G$ where here $cp(d)$ means the complexity of $d$ with respect to the fact that $D_U$ is 
the rational closure of $F[U,\eta,\alpha]$ and therefore of $F^\times X_U$ in $D$. This implies that 
the subring of $D$ which is generated by $D_U$ and $t$ is the set of all left (and right) polynomials 
in $t$ over $D_U$ which will be denoted by $D_U[t]$. We say that $t$ is {\it (left) algebraic} over 
$D_U$ if there exist $k\in\N$ and $d_0,d_1,\dots,d_{k-1}\in D_U$ such that 
$t^k+d_{k-1} t^{k-1}+\dots+d_1 t+d_0 = 0$ where {\it (right) algebraic} is defined accordingly. Of course 
both properties are equivalent for $t$ and if $t$ is not algebraic we call $t$ transcendental over 
$D_U$. In this case $D_U[t]$ is a skew polynomial ring in the indeterminate $t$ over $D_U$. 
\begin{definition} \label{Hughes-free}
Let the notation be as above. A division ring of fractions $D$ of a crossed product $R=F[G,\eta,\alpha]$ 
of a left-ordered group $G$ over a skew field $F$ is called Hughes-free if for any non-trivial finitely 
generated subgroup $H$ of $G$, any normal subgroup $N$ of $H$, and any $h \in H$ such that 
$H/N = \langle hN \rangle$ is an infinite cyclic group each element $t \in F^\times x_h$ is 
transcendental over $D_N$.
\end{definition}

{\bf Remarks.}
\begin{enumerate}
\item The definition of (Hughes-)free division rings of fractions in the original version appeared in 
      \cite{H1} for the first time where $G$ is assumed to be locally indicabel. 
\item If $R=F[G,\eta,\alpha]$ is given as above then it is not clear in general which basis $X_G$ of 
      the $F$-vector space $F[G,\eta,\alpha]$ is meant here. For example, in Proposition \ref{Einheitengruppe} 
      $R$ is a crossed product with respect to $X_G$ and $X'_G$. Moreover, in the definition above 
      besides the group $G$ also the basis $X_G$ occurs which leads to some kind of ambiguity. Later 
      on we will see that we have to be more careful in order to avoid misunderstandings and we have 
      to state precisely which basis is used.
\item It is easy to see that in Definition \ref{Hughes-free} the sentence {\it each element 
      $t \in F^\times x_h$ is transcendental over $D_N$} can be replaced by the sentence {\it there is 
      an element $t \in F^\times x_h$ which is transcendental over $D_N$}.
\item In \cite{JL}  it is shown that any group ring $F[G]$ of a locally indicable group $G$ over a field 
      $F$ with characteristic $0$ has a Hughes-free division ring of fractions. The question whether this 
      is true for an arbitrary crossed product $F[G,\eta,\alpha]$ of a locally indicable group $G$ is 
      still open.\\
\end{enumerate}
\begin{proposition}\label{allgemeinHughes-free}
Let $D$ be a division ring of fractions of a crossed product $F[G,\eta,\alpha]$ of a left-ordered group 
$G$ over a skew field $F$. Then the following statements are equivalent:
\begin{enumerate}
\item $D$ is Hughes-free.
\item If $H$ is a subgroup of $G$ and $N$ a normal subgroup of $H$ such that $H/N = \langle hN \rangle$ 
      is an infinite cyclic group then each $t \in F^\times x_h$ is transcendental over $D_N$.
\item If $H$ is a subgroup of $G$, $N$ a normal subgroup of $H$ such that $H/N = \langle hN \rangle$ 
      is an infinite cyclic group and if $h_1, \dots ,h_n \in H$ such that $h_1N, \dots ,h_nN$ are 
      pairwise distinct then $x_{h_1}, \dots , x_{h_n}$ are linearly independent over $D_N$.
\end{enumerate}
\end{proposition}

{\bf Proof.} $"(1) \Rightarrow (2):"$ We assume that there is an $a \in F^\times$ such that $t=ax_h$ 
is algebraic over $D_N$, that is, there exist $k\in \N$ and $d_0,d_1,\dots, d_{k-1} \in D_N$ 
satisfying
\[t^k + d_{k-1}t^{k-1} + \dots + d_1 t + d_0 =0.\]
Clearly, $d_0,\dots,d_{k-1}$ are elements of the rational closure of $F[U,\eta,\alpha]$ in $D$ where 
$U = \langle g_1,\dots,g_m \rangle$ is a suitable finitely generated subgroup of $N$. We consider 
the subgroup $H'$ of $H$ which is generated by $g_1,g_2,\dots,g_m$, and $h$. Then, $H'$ is a finitely 
generated subgroup of $G$ and $N':=H' \cap N$ is a normal subgroup of $H'$ such that 
$g_1,g_2,\dots,g_m \in N'$ and $h\not\in N'$. Therefore, $H'/N' = \langle hN' \rangle$ is an infinite 
cyclic group. By construction, $d_0,d_1,\dots,d_{k-1}$ belong to the rational closure $D_{N'}$ of 
$F[N',\eta,\alpha]$ in $D$. Since $D$ is a Hughes-free division ring of fractions of $F[G,\eta,\alpha]$ 
we get the desired contradiction.\\
$"(2) \Rightarrow (3):"$ Because of $H/N = \langle hN \rangle$ there exist $m_1, \dots ,m_n \in \Z$ 
and non-zero $a_1, \dots ,a_n \in D_N$ such that $x_{h_i} = a_i x^{m_i}_h$ for all $i=1,\dots,n$. 
Moreover, $m_1,\dots,m_n$ are pairwise distinct since $h_1N, \dots ,h_nN$ are pairwise distinct. Now, 
if $d_1,\dots,d_n \in D_N$ are given with $d_1x_{h_1} + \dots + d_nx_{h_n} = 0$ then 
$d_1a_1x^{m_1}_h + \dots + d_na_nx^{m_n}_h = 0$, that is $d_1 = \dots = d_n = 0$ for $x_h$ is 
transcendental over $D_N$ and $a_1, \dots ,a_n \not= 0$.\\
$"(3) \Rightarrow (1):"$ Let $H$, $N$, and $h$ are assumed as given in Definition \ref{Hughes-free} 
and let $d_0, \dots ,d_k$ are elements of $D_N$ with $d_kx^k_h + \dots + d_1x_h + d_0 = 0$. Since $hN$ 
has infinite order in $H/N$ the cosets $h^kN, \dots , hN, N$ are pairwise distinct. Thus, 
$x^k_h, \dots , x_h, x_e$ are linearly independent over $D_N$, that is $d_k = \dots = d_0 = 0$.\\
\qed\\

Proposition \ref{allgemeinHughes-free} shows that the two definitions of Hughes-freeness presented in 
\cite{H1} and \cite{L} for locally indicable groups coincide. It also gives rise to the following 
definition.

\begin{definition}\label{stronglyHughes-free}
Let $D$ be a division ring of fractions of a crossed product $F[G,\eta,\alpha]$ of a left-ordered 
group $G$ over a skew field $F$. Then $D$ is called strongly Hughes-free if the following holds true: 
If $H$ is an arbitrary subgroup of $G$ with normal subgroup $N$ of $H$ and if $h_1, \dots , h_n \in H$, 
$n \in \N$ are given such that $h_1N, \dots , h_nN$ are pairwise different, then $x_{h_1}, \dots , x_{h_n}$ 
are linearly independent over the rational closure $D_N$ of $F[N,\eta,\alpha]$ in D.
\end{definition}

The definition of strongly Hughes-free division rings of fractions in the original version appeared in 
\cite{L} where $G$ is assumed to be locally indicabel. The following proposition holds true more generally 
but in Section \ref{Hughes Theorems} we just need this version.

\begin{proposition}\label{freifuervarphi(R)}
Let the notation be as in Proposition \ref{Einheitengruppe} and let $D$ be a (strongly) Hughes-free 
division ring of fractions of $R=F[G,\eta,\alpha]$ with respect to $X_G$. Then $D$ is also (strongly) 
Hughes-free with respect to $X'_G$.
\end{proposition}

{\bf Proof.} We just consider the case where $D$ is strongly Hughes-free. For any $g \in G$ there are 
uniquely determined $a_g \in F^\times$ and $g' \in G$ such that $x'_g = \varphi(x_g) = a_gx_{g'}$. The 
mapping $\psi: G \longrightarrow G, g \longmapsto g'$ is bijective and for all $g,h \in G$ we get 
$\varphi(x_{gh}) = a_{gh}x_{\psi(gh)}$ and 
$\varphi(x_{gh})=\varphi(\eta(g,h)\inv) a_g \alpha_{\psi(g)}(a_h)\eta(\psi(g),\psi(h))x_{\psi(g)\psi(h)}$. 
Thus, $\psi$ is a group automorphism. In order to simplify the notation let an apostrophe indicate the 
image of a term under $\psi$. Now, let $H$, $N$, and $h_1, \dots ,h_n \in H$ are given as in 
Definition \ref{stronglyHughes-free}. Then $H'$ is a subgroup of $G$, $N'$ a normal subgroup of $H'$, and 
$h'_1 N',\dots,h'_n N'$ are pairwise distinct. Clearly, $F^\times X'_{N}$ and $F^\times X_{N'}$ generate 
the same subring of $D$ such that $D_{N'}$ is also the rational closure of $F^\times X'_{N}$ 
in $D$. By assumption, $a_{h_1}x_{h'_1},\dots,a_{h_n}x_{h'_n}$ are linearly independent over $D_{N'}$ where 
$a_{h_i}x_{h'_i} = x'_{h_i}$ for all $i=1,\dots,n$.\\
\qed\\

Next we introduce another notion of freeness in case that $\leq$ is a Conradian left-order of $G$. 
For any convex jump $(C',C)$ of $G$ with respect to $\leq$ let $D^+_C$ be the rational closure of 
$F[C,\eta,\alpha]$ and $D^-_C$ the rational closure of $F[C',\eta,\alpha]$ in $D$. If 
$g \in C \setminus C'$, that is $C=C^+_g$ and $C'=C^-_g$, we will also write $D^+_g$ and $D^-_g$ 
instead of $D^+_C$ and $D^-_C$ respectively. Then $tdt\inv \in D^-_g$ for all $d \in D^-_g$ and 
$t \in F^\times x_g$ since $C^-_g$ is normal in $C^+_g$.

\begin{definition}\label{free}
Let $D$ be a division ring of fractions of the crossed product $F[G,\eta,\alpha]$ of a group $G$ 
with Conradian left-order $\leq$ over a skew field $F$. Then $D$ is called free with respect to $\leq$ 
if for any convex jump $(C',C)$ of $G$ with respect to $\leq$ the following holds true: If 
$h_1, \dots , h_n \in C$, $n \in \N$ are given such that $h_1C', \dots , h_nC'$ are pairwise distinct, 
then $x_{h_1}, \dots , x_{h_n}$ are linearly independent over $D^-_C$.
\end{definition}

\begin{proposition}\label{hughes-free ist free}
Let $D$ be a division ring of fractions of the crossed product $R=F[G,\eta,\alpha]$ where $G$ is a 
locally indicable group and $F$ is a skew field. If $D$ is strongly Hughes-free then $D$ is free 
with respect to any Conradian left-order $\leq$ of $G$ and if $D$ is Hughes-free then $D$ is free 
with respect to any Conradian left-order $\leq$ of $G$ with maximal rank.
\end{proposition}

{\bf Proof.} If $D$ is strongly Hughes-free then the statement is obvious. Thus, let $\leq$ be a Conradian 
left-order of $G$ with maximal rank, let $(C',C)$ be a convex jump of $G$ with respect to $\leq$, and 
let $h_1, \dots , h_n \in C$, $n \in \N$ are given such that $h_1C', \dots , h_nC'$ are pairwise different. 
We define $H = H'C'$ where $H'$ is the subgroup of $C$ generated by $h_1,\dots,h_n$. Then $C'$ is a 
normal subgroup of $H$ and $H/C'$ is the subgroup of $C/C'$ generated by $h_1C', \dots , h_nC'$. Hence, 
$H/C'$ is cyclic since $C/C'$ is isomorphic to a subgroup of the additive group $\Q$. If $H=C'$ then 
$h_1,\dots ,h_n \in C'$, that is, $n=1$ and we are done. Otherwise, Proposition \ref{allgemeinHughes-free} 
can be applied.\\
\qed\\

The second part of Proposition \ref{hughes-free ist free} also appears in \cite{J}.\\

\section{A Structure Theorem for Free Division Rings of Fractions}\label{Struktursatz}

\vspace*{0.3cm}

Throughout this section we shall fix the following notation. $G$ is a (locally indicable) group with 
Conradian left-order $\leq$ and $F[G,\eta,\alpha]$ is a crossed product of $G$ over a skew field $F$ 
having a division ring of fractions $D$ which is free with respect to $\leq$. The convex subgroups of 
$G$ with lower neighbour in $\Co_\leq$ shall be called convex successors (of $G$ with respect to $\leq$). 
They form a subset $\Co^\ast_\leq$ of $\Co_\leq$ and if $C \in \Co^\ast_\leq$ is given then $C'$ always 
denotes the lower neighbour of $C$ such that $(C',C)$ is a convex jump. For example, $C_g^+$ is a convex 
successor whenever $g \not= e$. We intend to continue using the notation $D^+_C, D^-_C$, $D^+_g$, and 
$D^-_g$ as introduced before Definition \ref {free}. For any $C \in \Co^\ast_\leq$ we also 
fix a transversal $\frak C$ for $C'$ in $C$ satisfying $e \in \frak C$ and $\ofC$, $\tfC$ are defined 
analogously for $\oC, \tC \in \Co^\ast_\leq$. According to Proposition \ref{RGverschraenktueberRN} the 
ring $F[C,\eta,\alpha]$ now turns into a crossed product of $C/C'$ over $F[C',\eta,\alpha]$ with respect 
to the basis $X_{C/C'} = \{x_g \mid g \in {\frak C}\}$. Since $D$ is free with respect to $\leq$ 
the elements of $X_{C/C'}$ are linearly independent over $D^-_C$. This means that the subring of $D$ 
generated by $D^-_C$ and $X_{C/C'}$ is a crossed product of $C/C'$ over $D^-_C$ with respect to the basis 
$X_{C/C'}$ which will be denoted by $D^-_C[C/C',\eta,\alpha]$. Clearly, $D^+_C$ is the rational closure of 
$D^-_C[C/C',\eta,\alpha]$ in $D$. Because of Proposition \ref{abelschwohlgeordnet} and Proposition 
\ref{conradore} the ring $D^-_C[C/C',\eta,\alpha]$ is an Ore domain for $C/C'$ is abelian and therefore 
$D^+_C$ is the corresponding Ore division ring of fractions. Next we consider $D^-_C[C/C',\eta,\alpha]$ 
as a subring of the skew field $D^-_C((C/C',\eta,\alpha))$ of all formal power series in $C/C'$ over 
$D^-_C$ as given in \cite{Ne} that also contains $D^+_C$ as a subdivision ring in a unique manner because 
$D^-_C[C/C',\eta,\alpha]$ meets the Ore condition. Hence, any element $x \in D^+_C$ possesses a unique 
representation $x = \sum_{g \in \frak C} d_g x_g$ as a formal power series with $d_g \in D^-_C$ for all 
$g \in \frak C$ such that the support $\{g \in \frak C \mid d_g \not=0\}$ is well-ordered with respect 
to the given Conradian left-order of $G$.\smallskip\\
Unless otherwise specified, the complexity $cp(d)$ of an element $d \in D$ always refers to the way $d$ 
is built from elements of $\G = F^\times X_G$.

\begin{theorem}\label{DarstellungeninfreienKoerpern}
With the assumptions stated above, for any $x \in D$ with $cp(x) > 1$ the following hold true:
\begin{enumerate}
\item There exist $h \in G$ and $C \in C^\ast_\leq$ such that $x_h^{-1}x \in D^+_C$ and
      $x_h^{-1}x = \sum_{g \in \frak C} d_g x_g$ where $d_g \in D^-_C$ with $d_g \lhd x$ for all $g \in \frak C$.
\item There exist $h \in G$ and $C \in C^\ast_\leq$ such that $xx_h^{-1} \in D^+_C$ and 
      $xx_h^{-1} = \sum_{g \in \frak C} d_g x_g$ where $d_g \in D^-_C$ with $d_g \lhd x$ for all $g \in \frak C$.
\end{enumerate}
\end{theorem}

{\bf Remarks.}
\begin{enumerate}
\item We call $x_h^{-1}x = \sum_{g \in \frak C} d_g x_g$ as given above a left representation of $x$ and then  
      $x_{\oh h}^{-1}\gap (x_\oh x) = \sum_{g \in \frak C} (\alpha\inv_{\oh h}(\eta(\oh,h))d_g)\gap x_g$ 
      is a left representation of $x_\oh x$ for all $\oh \in G$. Similarly, 
      $xx_h^{-1} = \sum_{g \in \frak C} d_g x_g$ as given above is called a right representation of $x$ and 
      $(xx_\oh)\gap x_{h\oh}^{-1} = \sum_{g \in \frak C} (d_g\alpha_g(\eta(h,\oh)))x_g$ is a right 
      representation of $xx_\oh$ for all $\oh \in G$.
\item We will see later in Theorem \ref{EindeutigkeitderDarstellung} that the left and right representations 
      of an element $x \in D$, $cp(x) > 1$ is unique in a certain sense.
\item Because of $d_g \lhd x$ for all $g \in \frak C$ there are at least two non-zero summands in any 
      left and right representation of $x$.\\
\end{enumerate}

{\bf Proof of Theorem \ref{DarstellungeninfreienKoerpern}.} We will prove both statements of the theorem 
simultaniously by transfinite induction on the complexity of $x$ and consider the case 
$x \in F[G,\eta,\alpha]$ first. Then
\[x = a_1x_{g_1} + a_2x_{g_2} + \dots + a_nx_{g_n}\]
where $a_1,\dots,a_n \in F^\times$ and $g_1, \dots , g_n \in G$ are pairwise distinct. Moreover, $n \geq 2$ 
because of $cp(x) >1$. We may choose notation so that $C^+_{g_1} \subseteq \dots \subseteq C^+_{g_n}$ 
and two different cases occur.\smallskip\\
Case 1: $C^+_{g_1} \subset C^+_{g_n}$. Then 
\[C^+_{g_1} \subseteq \dots \subseteq C^+_{g_m} \subset C^+_{g_{m+1}} = \dots = C^+_{g_n}\]
for some $m<n$. We put $C := C^+_{g_n}$, that is $C' := C^-_{g_n}$, and 
$d_e = a_1x_{g_1} + \dots + a_mx_{g_m} \in D^-_C$ with $d_e \lhd x$ since $m<n$. Collecting all summands 
in $a_{m+1}x_{g_{m+1}} + \dots + a_nx_{g_n}$ that correspond to the same cosets of $C'$ there exist 
pairwise distinct $h_1,\dots,h_r \in \frak C$ such that 
$a_{m+1}x_{g_{m+1}} + \dots + a_nx_{g_n} = d_{h_1}x_{h_1} + \dots + d_{h_r}x_{h_r}$ where 
$d_{h_i} \in D^-_C$ and $d_{h_i} \lhd x$ for all $i=1,\dots,r$. Finally, we put $h=e$ for both parts of 
the theorem and obtain
\[x_h^{-1}x = d_ex_e + d_{h_1}x_{h_1} + \dots + d_{h_r}x_{h_r} = xx_h\inv\]
as desired.\smallskip\\
Case 2: $C^+_{g_1} = C^+_{g_2} = \dots = C^+_{g_n}$. Because of $n\geq 2$ there exists at least 
one index $i$ such that $g_i$ is not trivial which yields that $g_i$ is not trivial for all $i=1,2,\dots,n$. 
Defining $h=g_1$ provides
\begin{center}
$x_h^{-1}\gap x = a'_1 x_e + a'_2x_{g^{-1}_1g_2} + \dots + a'_n x_{g^{-1}_1g_n}$,\\
$x\gap x_h^{-1} = \ta_1 x_e + \ta_2 x_{g_2 g^{-1}_1} + \dots + \ta_n x_{g_n g^{-1}_1}$
\end{center}
where $g^{-1}_1 g_i, g_ig^{-1}_1 \not= e$ for $i\geq 2$ and $a'_1,\dots,a'_n, \ta_1,\dots,\ta_n \in F^\times$. 
We can now proceed analogously to Case 1 again.\bigskip\\
We next turn to the case $x\not\in D_{[0]}$ and assume that the theorem has been proved for all elements 
in $D$ which are strictly simpler than $x$. First of all, let $x$ be an atom. Then $x^{-1} \lhd x$ by 
\cite[Proposition 4.1]{DGH}. There exist $h\in G$ and $C \in \Co^\ast_\leq$ such that 
$x_h^{-1}x^{-1} \in D^+_C$ and 
\[x_h^{-1}x^{-1} = \sum_{g \in \frak C} d_g x_g \mbox{\ where\ } d_g \in D^-_C \mbox{\  with\ } d_g \lhd x\inv 
                   \mbox{\ for all\ } g \in \frak C.\]
From this we will deduce a right representation $x x_h = \sum_{\og \in \frak C}\od_\og x_\og \in D^+_C$ 
of $x x_h$ where $\od_\og \in D^-_C$ and $\od_\og \lhd x$ for all $\og \in \frak C$ such that 
$x x\inv_{h\inv} = \sum_{\og \in \frak C}(\od_\og \alpha_\og(\eta(h\inv,h)\inv ))x_\og$ is a right representation 
of $x$ since $\od_\og \alpha_\og(\eta(h\inv,h)\inv ) \lhd x$ by \cite[Proposition 4.8]{DGH}. \smallskip\\
Now, let $g_0$ be minimal in the support of 
$x\inv_h x\inv$. For any $g$ from this support, $g > g_0$ there exist uniquely determined $g' \in \frak C$ 
and $a_g \in F^\times X_{C'}$ satisfying $g g\inv_0 C' = g'C'$ and $x_g x\inv_{g_0} = a_{g}x_{g'}$. 
Furthermore, $g'C' = g\inv_0gC'$ and $g\inv_0g > e$ for $C/C'$ is abelian, that is $g' > e$ and 
$\{g' \mid d_g \not = 0\}$ is well-ordered with respect to the given order $\leq$ of $G$. Therefore,
\[1+\sum_{g' > e} d\inv_{g_0} d'_{g'}x_{g'} \in D^-_C((C/C',\eta,\alpha))\]
where $d'_{g'} := d_g a_{g}$ and
\begin{center}
$x_h^{-1} \gap x^{-1} = d_{g_0}( 1+\sum_{g' > e} d\inv_{g_0} d'_{g'}x_{g'})x_{g_0}$
\end{center}
such that
\begin{center}
$x \gap x_h = x\inv_{g_0}( 1+\sum_{g' > e} d\inv_{g_0} d'_{g'}x_{g'})\inv  d^{-1}_{g_0}$.
\end{center}
According to \cite[Section4]{Ne} there exists a power series representation
\[( 1+\sum_{g' > e} d\inv_{g_0} d'_{g'}x_{g'})\inv = 1 + \sum_{\tg>e} \td_\tg x_\tg,\]
where $\td_\tg \in D^-_C$ for any $\tg>e$ and next we verify $\td_\tg \lhd x$. The proof of 
\cite[Theorem 4.9]{Ne} shows that any $\td_\tg$ is a sum of products with factors from
\[\Lambda = \{\pm 1, d\inv_{g_0}, d'_{g'}, x _{g'}, x\inv_{g'} \mid g' \in \frak C\} \cup F^\times X_{C'}.\]
Clearly, $\pm 1$, all $x_{g'}, x\inv_{g'},$ and all elements from $F^\times X_{C'}$ are strictly simpler 
than $x$. Moreover, \cite[Proposition 4.1]{DGH} and $d_g \lhd x^{-1}$ imply $d_g, d^{-1}_{g} \lhd x$ 
such that $d^{-1}_{g_0} \lhd x$ and $d'_{g'} = d_g a_{g}\lhd x$ by \cite[Proposition 4.8]{DGH}. 
Because of \cite[Proposition 4.1]{DGH} any product of elements from $\Lambda$ is strictly simpler than $x$ 
which yields $\td_\tg \lhd x$.
This means
\[x \gap x_h = x\inv_{g_0}(1 + \sum_{\tg>e} \td_\tg x_\tg)d^{-1}_{g_0}
             = x\inv_{g_0}d^{-1}_{g_0} + \sum_{\tg>e} x\inv_{g_0}\td_\tg x_\tg d^{-1}_{g_0}.\]
This presentation of $x \gap x_h$ can now easily be transformed into a formal power series of the desired 
form since $x\inv_{g_0}d^{-1}_{g_0} x_{g_0}$ and $x\inv_{g_0}\td_\tg x_\tg d^{-1}_{g_0}x\inv_\tg x_{g_0}$ 
are in $D^-_C$ for all $\tg>e$ and they are all strictly simpler than $x$. Similar arguments show that $x$ 
also possesses a left representation if $x$ is an atom.\bigskip\\
We continue with the case that $x$ is additively indecomposable and consider the complete multiplicative 
decomposition $x=x_1 \cdot x_2 \cdot \dots \cdot x_n$ of $x$ with $n \in \N$. According to 
\cite[Theorem 4.6]{DGH} any $x_i$, $i=1,\dots,n$ is a proper atom. We prove the claim of the theorem 
for $x$ by induction on $n$ using that the theorem holds true for all elements in $D$ which are strictly 
simpler than $x$. If $n=1$ then $x$ is a proper atom which has been treated above. Thus, let $n>1$. 
We define $\ox := x_2 \cdot \dots \cdot x_n$ which is additively indecomposable by \cite[Lemma 4.5]{DGH} 
and $\ox = x_2 \cdot \dots \cdot x_n$ is a complete multiplicative decomposition of $\ox$ by 
\cite[Theorem 4.6]{DGH}. We apply the induction hypothesis to obtain $\oh \in G$, $\oC \in \Co^\ast_\leq$ 
such that $x_\oh^{-1}\gap \ox \in D^+_\oC$ and
\begin{center}
$x_\oh^{-1}\gap \ox = \sum_{\og \in \ofC} \od_\og x_\og$\ with\ $\od_\og \in D^-_\oC$ \ and\ 
$\od_\og \lhd \ox$\ for all $\og \in \ofC$.
\end{center}
Since $x_1 x_\oh$ is a proper atom there exist $\tth \in G$, $\tC \in \Co^\ast_\leq$ such that 
$x_\tth^{-1}\gap (x_1\gap x_\oh) \in D^+_\tC$ and
\begin{center}
$x_\tth^{-1}\gap (x_1\gap x_\oh) = \sum_{\tg \in \tfC} \td_\tg x_\tg$\ with\ $\td_\tg \in D^-_\tC$ \ and\ 
$\td_\tg \lhd x_1\gap x_\oh$\ for all $\tg \in \tfC$.
\end{center}

We shall distinguish three cases.\smallskip\\
Case 1: $\tC \subset \oC$, that is $\tC \subseteq \oC'$, and $D^+_\tC$ is a subdivision ring of $D^-_\oC$. 
Further,
\[x_\tth^{-1}\gap x = x_\tth^{-1} x_1 \ox = x_\tth^{-1} x_1 x_\oh \gap x_\oh^{-1} \ox = 
    \sum_{\og \in \ofC} (x_\tth^{-1} x_1 x_\oh \gap \od_\og) x_\og \mbox{\ with\ } 
    x_\tth^{-1} x_1 x_\oh \gap \od_\og \in D^-_\oC.\]
Since $x_\tth^{-1} \gap x_1 \gap x_\oh$ and $x_1$ have the same complexity we conclude 
$x_\tth^{-1}\gap x_1 x_\oh \gap \od_\og \lhd x_1 \ox = x$ because of 
$\od_\og \lhd \ox$ and \cite[Theorem 4.6]{DGH}.\smallskip\\
Case 2: $\oC \subset \tC$, that is $\oC \subseteq \tC'$, and $D^+_\oC$ is a subdivision ring of $D^-_\tC$. 
From $x_\oh^{-1} \ox \in D^-_\tC$ we get $x_\tg x\inv_\oh \ox x\inv_\tg \in D^-_\tC$ 
and $\td_\tg x_\tg x\inv_\oh \ox x\inv_\tg \in D^-_\tC$ for all $\tg \in \tfC$ where, in addition   
$cp(x_\tg x_\oh^{-1} \gap \ox \gap x\inv_\tg) = cp(\ox)$ and $cp(\td_\tg) < cp(x_1 x_\oh) = cp(x_1)$ 
imply $\td_\tg \gap x_\tg \gap x_\oh^{-1} \gap \ox \gap x\inv_\tg \lhd x_1 \ox = x$ due to 
\cite[Theorem 4.6]{DGH}. Thus,
\[x_\tth^{-1} \gap x = x_\tth^{-1} \gap x_1 \ox = 
    x_\tth^{-1} \gap x_1 \gap x_\oh \gap x_\oh^{-1} \gap \ox = 
    \sum_{\tg \in \tfC} (\td_\tg x_\tg x\inv_\oh \ox x\inv_\tg) x_\tg\]
is a left representation of $x$.\smallskip\\
Case 3: $\oC = \tC$. We define $C:= \oC = \tC$ and obtain 
$x\inv_\oh \ox, x\inv_\tth (x_1 x_\oh) \in D^+_C$ with left representations 
$x_\oh^{-1}\gap \ox = \sum_{g \in \frak C} \od_g x_g$ and 
$x_\tth^{-1}\gap (x_1\gap x_\oh) = \sum_{g \in \frak C} \td_g x_g$ where $\od_g, \td_g \in D^-_C$ such 
that $\od_g \lhd \ox$, $\td_g \lhd x_1 x_\oh$. In the skew field $D^-_C((C/C',\eta,\alpha))$ of all 
formal power series in $C/C'$ over $D^-_C$ we obtain
\[x\inv_\tth x = x\inv_\tth x_1 x_\oh x\inv_\oh \ox = (\sum_{g \in \frak C} \td_g x_g)
                                                      (\sum_{g \in \frak C} \od_g x_g)\]
as a product of two formal power series. The resulting series has coefficients which are sums of products 
with two factors each where the first factor always is simpler than $x_1$ and the second is simpler 
than $\ox$. According to \cite[Theorem 4.6]{DGH} all these products are simpler than $x_1\ox = x$. Since $x$ 
is additively indecomposable any sum of these products is simpler than $x$ by \cite[Proposition 3.1]{DGH}. 
This provides the desired left representation of $x = x_1 \ox$. Similar arguments can be used to show that 
$x$ has a right representation as stated in the theorem.\bigskip\\
We finally turn to the general case where $x$ possesses a complete additive decomposition 
$x = x_1 + x_2 + \dots + x_n$ and we prove the claim of the theorem for $x$ by induction on $n$. If $n=1$ 
then $x$ is additively indecomposable which has been treated above. Thus, let $n>1$. 
\cite[Theorem 3.6]{DGH} ensures that each $x_i$ is additively indecomposable and 
$\ox := x_2 + \dots + x_n$ is a complete additive decomposition of $\ox$ by \cite[Theorem 3.6]{DGH}. 
We will admit for a moment that $cp(x_1), cp(\ox) > 1$ which enable us to apply the induction hypothesis. 
Hence, there exist $\oh \in G$, $\oC \in \Co^\ast_\leq$ such that $x_\oh^{-1} \gap \ox \in D^+_\oC$ and
\[x_\oh^{-1} \gap \ox = \sum_{\og \in \ofC} \od_\og \gap x_\og \mbox{\ with\ } \od_\og \in D^-_\oC 
 \mbox{\ and\ } \od_\og \lhd \ox \mbox{\ for all\ } \og \in \ofC.\]

In addition, there are $\tth \in G$, $\tC \in \Co^\ast_\leq$ satisfying $x_\tth^{-1} \gap x_1 \in D^+_\tC$ 
and
\[x_\tth^{-1} \gap x_1 = \sum_{\tg \in \tfC} \td_\tg \gap x_\tg \mbox{\ with\ } 
 \td_\tg \in D^-_\tC \mbox{\ and\ } \td_\tg \lhd x_1 \mbox{\ for all\ } \tg \in \tfC.\]
We just investigate the case $\tC \subseteq \oC$. If $\oC \subseteq \tC$ we proceed 
similarly. Again, several subcases have to be distinguished.\smallskip\\
Case 1: $\tC = \oC$. We put $h := \oh\inv \tth$ and consider two subcases where the first one 
refers to $h \in \tC = \oC$. Then, we define $C := \tC = \oC$ and write
\[x_\oh^{-1} \gap \ox = \sum_{g \in \fC} \od_g \gap x_g,\ 
 x_\tth^{-1} \gap x_1 = \sum_{g \in \fC} \td_g \gap x_g \mbox{\ with\ } 
 \od_g, \td_g \in D^-_C \mbox{\ and\ } \od_g \lhd \ox, \td_g \lhd x_1.\]
For all $g \in \fC$ there exist uniquely determined $g' \in \fC$ and $u_g \in F^\times X_{C'}$ such 
that $x_hx_g = u_g x_{g'}$ and $\{g' \mid \td_g \not= 0\}$ is well-ordered with respect to $\leq$. 
We now introduce $\td '_{g'} := (x\inv_\oh (\eta(\oh,h)\inv)x_\oh)( x_h \td_g x\inv_h) u_g \in D^-_C$ 
for all $g \in \fC$ where $\td_g$ and $\td'_{g'}$ have the same complexity. A straightforward calculation 
shows
\[x_\oh^{-1} \gap x_1 = \sum_{g' \in \fC} \td'_{g'} \gap x_{g'} \mbox{\ with\ } 
 \td'_{g'} \in D^-_C \mbox{\ and\ } \td'_{g'} \lhd x_1.\]
Therefore, in this subcase we may assume $\tth = \oh$. Adding the power series representations 
of $x_\oh^{-1} \gap \ox$ and $x_\oh^{-1} \gap x_1$ yields a representation of $x\inv_\oh (x_1 + \ox)$ 
as a power series $\sum_{g \in \fC} d_g x_g$ with coefficients $d_g = \td_g + \od_g \in D^-_C$ where 
$\td_g \lhd x_1$ and $\od_g \lhd \ox$, that is, $d_g \lhd x$ because of \cite[Theorem 3.6]{DGH}.
\smallskip\\
Now, we turn to the second subcase of $\tC=\oC$ that refers to $h \not\in \oC$ and define 
$C=C^+_h \in \Co^\ast_\leq$. This provides $\tC, \oC \subseteq C'$ which means that $D^+_\tC$ and 
$D^+_\oC$ are subdivision rings of $D^-_C$. Thus, $x\inv_\oh \ox, x\inv_\tth x_1 \in D^-_C$ follows. 
Since $h \in C\setminus C'$ there exist $g \in \fC$, $g \not= e$ and $u \in F^\times X_{C'}$ satisfying 
$x\inv_h = u x_g$ such that 
\[x\inv_\tth x = x\inv_\tth x_1 + x\inv_\tth \ox = d_e x_e + d_g x_g\]
with $d_e = x\inv_\tth x_1$ and $d_g = (x\inv_h((x\inv_\oh \eta(\oh,h)x_\oh)(x\inv_\oh \ox))x_h)u$. 
Because of $d_e, d_g \in D^-_C$, $cp(d_e) = cp(x_1) < cp(x)$, and $cp(d_g) = cp(\og) < x$ this is a 
left representation of $x$.\smallskip\\
Case 2: $\tC \subset \oC$. Again, let $h=\oh\inv \tth$. As above we investigate the first subcase 
where $h \in \oC$ and choose $C := \oC$, that is, $\tC \subseteq C'$. Then, 
$x\inv_\tth x_1 \in D^+_\tC \subseteq D^-_C$. Because of $h \in C$ there exist $g' \in \fC$ and 
$u \in F^\times X_{C'}$ such that $x_h = u x_{g'}$. Now, we define 
$d'_{g'} = ((x\inv_\oh \eta(\oh,h)\inv x_\oh)(x_h (x\inv_\tth x_1) x\inv_h))u$ and conclude 
$d'_{g'} \in D^-_C$ where in addition $cp(d'_{g'}) = cp(x_1)$. A straightforward calculation shows 
$x\inv_\oh x_1 = d'_{g'} x_{g'}$ such that
\[x\inv_\oh x = x\inv_\oh x_1 + x\inv_\oh \ox = d'_{g'} x_{g'} + \sum_{g \in \fC} \od_g x_g = 
\sum_{g \in \fC} d_g x_g\]
where $d_g \in \{0, \od_g, d'_{g'}, d'_{g'} + \od_g\} \in D^-_C$ for all $g \in \fC$ and 
$d_g \lhd x$ by \cite[Theorem 3.6]{DGH}.\smallskip\\
The second subcase of $\tC \subset \oC$ deals with $h \not\in \oC$ and we proceed exactly as in the 
second subcase of Case 1.\smallskip\\
Recall, that $cp(x_1), cp(\ox) > 1$ has been presumed at the intermediate stage of the proof where 
$x$ is assumed to be a general element of $D$ having a proper additive decomposition. Clearly, the 
case $cp(x_1) = cp(\ox) = 1$ has already been treated when the theorem has been shown for 
$x \in F[G,\eta,\alpha]$. Thus, let us suppose for example that $cp(x_1) = 1$ and $cp(\ox) > 1$. 
Then, $x_1 = ax_g$ with $a \in F^\times$ and $g \in G$ which means $x\inv_\tth x_1 = d_e x_e$ where 
$\tth = g$ and $d_e = x\inv _g a x_g \in F^\times$. Even though $x\inv_\tth x_1 = d_e x_e$ is not a 
proper left representation we introduce formally $\tC = \{e\}$ such that $\tC \subset \oC$ since 
$cp(\ox) > 1$ and proceed as in Case 2.\smallskip\\
Similar arguments yield that $x$ also has a right representation.\\
\qed\\
\begin{theorem}\label{EindeutigkeitderDarstellung}
With the assumptions as above, let $x \in D$ with $cp(x) > 1$ such that 
$x_\oh^{-1} \gap x \in D^+_\oC$ for some $\oh \in G$, $\oC \in \Co^\ast_\leq$ and with respect to this 
let $x_\oh^{-1} \gap x = \sum_{\og \in\ofC} \od_\og \gap x_\og$ such that $\od_\og \in D^-_\oC$ 
for all $\og \in \ofC$ and $\od_\og \not=0$ for at least two different $\og \in \ofC$. Then 
the following hold true where the notation is as in the first part of Theorem 
\ref{DarstellungeninfreienKoerpern}:
\begin{enumerate}
\item $C = \oC$.
\item $hC = \oh C$.
\item If $c \in C$ with $h = \oh c$  and if $u_g \in F^\times X_{C'}$, $g' \in \fC$ with 
$x_cx_g = u_g x_{g'}$ for all $g \in \fC$ then 
$\od_{g'} = \alpha\inv_\oh(\eta(\oh,c)\inv) x_c d_g x\inv_c u_g$ for all $g \in\fC$.
\item $\od_g \lhd x$ for all $g \in \fC$ and 
      $x_\oh^{-1} \gap x = \sum_{g \in \fC} \od_g \gap x_g$ is a left representation of $x$.
\end{enumerate}
\end{theorem}

{\bf Proof.} (1): We begin with the special case $\oh = e$ and assume $C \subset \oC$. Then 
$C \subseteq \oC'$ and therefore $x\inv_h x \in D^+_C \subseteq D^-_\oC$.\smallskip\\
Case 1: $C^+_h \subset \oC$. Then $C^+_h \subseteq \oC'$ and $x\inv_\oh x = x = x_h(x_h^{-1}x) \in D^-_\oC$. 
On the other hand $x_\oh^{-1} \gap x = \sum_{\og \in\ofC} \od_\og \gap x_\og \in 
D^-_\oC((\oC/\oC', \eta, \alpha))$ such that $\od_\og \not=0$ for at least two different $\og \in \ofC$ 
which is a contradiction.\smallskip\\
Case 2: $C^+_h = \oC$, that is $C^-_h = \oC'$, and $x x\inv_h = x_h(x\inv_h x) x\inv_h \in D^-_\oC$ 
follows. There exist $h' \in \ofC$, $u_h \in F^\times X_{\oC'}$ satisfying $x_h =u_h x_{h'}$. 
Hence, $x\inv_\oh x = x = (x x\inv_h u_h) x_{h'}$ with $x x\inv_h u_h \in D^-_\oC$ is a power 
series representation of $x\inv_\oh x$ which is a contradiction as in Case 1.\smallskip\\
Case 3: $\oC \subset C^+_h$. Then $\oC \subseteq C^-_h$ follows and therefore 
$x_h^{-1}x \in D^+_C \subseteq D^-_h$. Moreover, $x = x\inv_\oh x \in D^+_\oC \subseteq D^-_h$ which implies 
$x\inv_h = (x\inv_h x)x\inv \in D^-_h$. Hence, $x_e$ and $x_h$ are linearly dependent over $D^-_h$, a 
contradition.\smallskip\\
Similar arguments disprove $\oC \subset C$ where it should be noted that $d_g \lhd x$ for 
$g \in \fC$ will not be used directly. Instead of this we apply the derived property $d_g \not=0$ for 
at least two different $g \in \fC$. All this shows the first statement of the theorem for the case 
$\oh = e$. Now, the general situation can be treated as follows. We put $x' = x\inv_\oh x$ 
with $x\inv_e x' = \sum_{\og \in\ofC} \od_\og \gap x_\og$ and obtain from 
$x\inv_h x = \sum_{g \in \fC} d_gx_g$ a left representation $x\inv_{h'} x' = \sum_{g \in \fC} d'_g x_g$ 
where $h' = \oh\inv h$ and $d'_g =  \alpha\inv_h(\eta(\oh,h')\inv)d_g$ for all $g \in \fC$.\smallskip\\
(2): Again, we first suppose $\oh = e$. Then Case 3 above shows that $C = \oC \subset C^+_h$ cannot hold 
true. Thus $h \in C$ follows which means $hC = \oh C = C$. In the general situation we also define 
$x' = x\inv_\oh x$ and consider $x\inv_e x'$, $x\inv_{h'} x'$. This provides $\oh\inv h = h' \in C$.
\smallskip\\
(3): Here we just have to transform the left representation $x\inv_h x = \sum_{g \in \fC} d_g x_g$ by 
means of $x\inv_h = x\inv_c \alpha\inv_\oh(\eta(\oh,c)) x\inv_\oh$ into a formal power series 
representation of $x\inv_\oh x$.\smallskip\\
(4): This is a direct consequence of statement (3) and \cite[Proposition 4.8]{DGH} since $d_g \lhd x$ 
for all $g \in \fC$.\\
\qed\\

{\bf Remarks.}
\begin{enumerate}
\item Similar arguments can be used in order to prove a corresponding theorem for right representations.
\item With the assumptions as above, let $H$ be a convex subgroup of $G$ and let $\leq$ also denote 
      the Conradian left-order of $H$ induced by the left-order $\leq$ of $G$. Then, the rational 
      closure $D_H$ of $F[H,\eta,\alpha]$ in $D$ is free with respect to $\leq$ since the convex 
      subgroups of $H$ coincide with the convex subgroups of $G$ contained in $H$. We will apply our 
      results to $D_H$ instead of $D$, that is, we replace $G$ by $H$. This means that for any $x \in D_H$, 
      $x \not\in F X_H$ there exist an $h \in H$, a convex subgroup $\oC$ of $H$ which is a convex successor 
      in $H$, and a left representation $x\inv_\oh x = \sum_{\og \in \ofC} \od_\og x_\og$ of $x$ where any 
      $\od_\og, \og \in \ofC$ is strictly simpler than $x$ where here simpler refers to the complexity 
      which describes how an element of $D_H$ is built from elements in $F[H,\eta,\alpha]$. Thus, 
      $\od_\og \not=0$ holds true for at least two different $\og \in \ofC$. We now apply Theorem 
      \ref{EindeutigkeitderDarstellung} and see that this left representation of $x$ as an element 
      from $D_H$ is also a left representation of $x$ as an element from $D$ since $\oC$ is a convex 
      successor in $G$ as well. Thus, any $\od_\og, \og \in \ofC$ is also strictly 
      simpler than $x$ where now simpler refers to the complexity originally defined for the elements 
      in $D$.\smallskip\\
\end{enumerate}

\begin{definition}\label{Definitionreduziert}
Let $x \in D$ such that $cp(x) > 1$ and let $x_h^{-1} \gap x = \sum_{g \in \fC} d_g x_g$ be as in 
Theorem \ref{DarstellungeninfreienKoerpern}. Then, $x$ is called reduced if $hC = C$, that is 
$h \in C$.
\end{definition}

Because of Theorem \ref{EindeutigkeitderDarstellung} the definition of {\it reduced} does not depend 
on the given left re\-presentation of $x$. If $x \in D$ is reduced then there is always a uniquely 
determined left representation $x = \sum_{g \in \fC} d_g x_g$ which we will call the {\it pure left 
representation} of $x$ and which is also a right representation. Furthermore, for any $x \in D$ with 
$cp(x)>1$ there exist an $h \in G$ and a reduced $x'\in D$ such that $x=x_h x'$.\\

{\bf Remarks about left representations of sums and products.}\smallskip\\
The following remarks apply in a similar manner to right representations. If $x,y \in D$ are given then 
left representations of $x+y$ and $x \cdot y$ essentially result from adding and multiplying left 
representations of $x$ and $y$ respectively as described in the proof of Theorem 
\ref{DarstellungeninfreienKoerpern}. We shall give a more detailed explanation of this which 
will be needed in Section \ref{Quotientenschiefkoerper und Potenzreihen}. To improve the readability 
we shall introduce left and right representations also for elements from $FX_G$. This will be done in 
an obvious way, that is, for $ax_g \in F X_G$ we call $ax_g$ itself a left and right representation 
and one shall always bear in mind that this does not have the proper meaning.\smallskip\\

{\underline {Left representations of sums}.}\smallskip\\
If $x$ and $y$ are in $F X_G$ then a left representation of $x + y$ results obviously. Thus, let us 
assume $cp(x) > 1$ or $cp(y) > 1$. In the proof of Theorem \ref{DarstellungeninfreienKoerpern} we have 
explained how to derive a left representation of $x = x_1 + \ox$ from those of $x$ and $\ox$ and this 
method can almost completely be applied to $x$, $y$, and $x+y$. The first observation is that two 
essentially different cases have to be distinguished where the second case has two subcases:\smallskip\\
$\bf A1$: Up to exchanging the roles of $x$ and $y$ there exist $h \in G$, $C \in \Co^\ast_\leq$, and 
$g \in \fC, g\not=e$ such that $x_h^{-1} \gap x, x_h^{-1} y x_g^{-1} \in D^-_C$ and
\[x_h^{-1} \gap (x+y) = x_h^{-1} \gap x + (x_h^{-1} y x_g^{-1}) \gap x_g\]
is a left representation of $x+y$.\smallskip\\
$\bf A2$: Up to exchanging the roles of $x$ and $y$ there exist $h \in G$ and $C \in \Co^\ast_\leq$ 
such that
\[x_h^{-1} \gap x = \sum_{g \in \fC} d_g x_g \mbox{\ and\ } 
  x_h^{-1} \gap y = \sum_{g \in \fC} d'_g x_g,\]
where the left equation shows a proper left representation of $x$. If the power series given on the 
right side has more than one non-trivial summand then it is a proper left representation of $y$. However, 
it is also possible that only one non-zero summand exists. But in any case, $d_g \lhd x$ and $d'_g \unlhd y$ 
hold true for all $g \in \fC$ and two subcases occur:\smallskip\\
$\bf A2.1$ There are at least two different $g \in \fC$ with $d_g + d'_g \not = 0$. Then, 
Theorem \ref{EindeutigkeitderDarstellung} ensures that
\[x_h^{-1} \gap (x+y) = \sum_{g \in \fC} (d_g + d'_g) \gap x_g\]
is a proper left representation of $x+y$.\smallskip\\
$\bf A2.2$ There is only one or no $g \in \fC$ with non-zero $d_g + d'_g$. Then,
\[x_h^{-1} \gap (x+y) = (d_g + d'_g) \gap x_g \mbox{\ or\ } x + y = 0.\]
In this case we do not obtain a proper left representation of $x+y$. However, in further investigations 
we will frequently use transfinite induction on the complexity of the elements of $D$ and then we will 
be able to apply even in this subcase the induction hypothesis to $x + y$.\\

{\underline {Left representations of products}.}\smallskip\\
For later use it will be enough to restrict our discussion to the case where $x$ and $y$ are reduced 
and $cp(x),cp(y) > 1$. In the proof of Theorem \ref{DarstellungeninfreienKoerpern} we have seen how 
to deduce a left representation of $x = x_1 \cdot \ox$ from those of $x$ and $\ox$ and again, these 
ideas can almost completely be transferred to $x$, $y$, and $x \cdot y$. Some arguments are actually 
becoming simpler since $x$ and $y$ are assumed to be reduced. As in the proof of 
Theorem \ref{DarstellungeninfreienKoerpern} there are three essentially different cases to be 
distinguished:\smallskip\\
$\bf M1$: There is a left representation $y = \sum_{g \in \fC} d_g \gap x_g$ of $y$ such that 
$x \in D^-_C$. Then, 
\[xy = \sum_{g \in \fC} (xd_g) \gap x_g\]
is a left representation of $xy$.\smallskip\\
$\bf M2$: There is a left representation $x = \sum_{g \in \fC} d_g \gap x_g$ of $x$ such that 
$y \in D^-_C$. Then,
\[xy = \sum_{g \in \fC} (d_g \gap x_g y x\inv_g) \gap x_g\]
is a left representation of $xy$.\smallskip\\
$\bf M3$: There are left representations  $x = \sum_{g \in \fC} \od_g \gap x_g$ and 
$y = \sum_{g \in \fC} \td_g \gap x_g$ of $x$ and $y$ respectively. Then,
\[xy = \sum_{g \in \fC} d_g \gap x_g \in D^-_C((C/C',\eta,\alpha)),\]
where
\[d_g x_g = \sum_{\og, \tg \in \fC \atop \og \tg C' = gC'} 
            \od_\og \gap x_\og \gap \td_\tg \gap x_\tg \mbox{\ \ and \ }
      d_g = \sum_{\og, \tg \in \fC \atop \og \tg C' = gC'} 
            \od_\og \big( (x_\og \td_\tg x\inv_\og) x_\og x_\tg x\inv_g\big)\]
for all $g \in \fC$ such that each $d_g$ is a finite sum with summands which are products with two factors 
each, where the left factor is strictly simpler than $x$ and the right factor strictly simpler than $y$. 
If $d_g$ is non-zero for at least two $g \in \fC$ then 
$xy = \sum_{g \in \fC} d_g \gap x_g$ is a left representation of $xy$. Otherwise, 
there exists a $g \in \fC$ such that $xy = d_g x_g$ with $d_g$ as above. This is no proper left 
representation in general but also this presentation will be needed in later applications.\\

\section{The Vector Space of Formal Power Series}\label{Formale Potenzreihen}

\vspace*{0.4cm}

In this section we summarize all propositions and theorems which refer to Dubrovin's method of embedding 
left-ordered groups into division rings and which will be used in this paper. They can be found for example 
in \cite{D1,D2,D3} where, however, some proofs are kept rather short and also results we need are 
occasionally hidden or they occur as special cases in a broader context. The author of this paper 
is preparing a book which is devoted to Dubrovin's work on this subject. It will provide a comprehensive 
representation of his results about the embedding of left-ordered groups into division rings by means 
of formal power series where emphasis is placed on detailed proofs. This involves for example all 
propositions and theorems included in this section. The notation we use here is similar to that in 
\cite{GS} and is also adapted to the subject of this paper. Therefore, it differ from what is introduced 
in \cite{D1,D2,D3}.   
\smallskip\\
Let $R=F[G,\eta,\alpha]$ be a crossed product of a group $G$ over a skew field $F$ where $G$ is 
endowed with a left-order $\leq$ which need not be Conradian. We also fix a basis $X_G=\{x_g\ |\ g \in G\}$ 
of the $F$-vector space $F[G,\eta,\alpha]$. A mapping
\[m: G \longrightarrow F, g \longmapsto m_g\]
is called a formal power series in $G$ over $F$ if its support $\supp m = \{g \in G \ |\ m_g \not=0\}$ 
is well-ordered with respect to $\leq$. We write $m=\Null$ in case of $m(g) = m_g = 0$ for all $g \in G$. 
The set of all formal power series in $G$ over $F$ will be denoted by $F((G))$.

\begin{definition}
If $m \in F((G))$, $m \not=\Null$ then $v(m)$ denotes the minimal element of $\supp m$.
\end{definition}

For $m = \Null$ we put $v(m) = \infty$ where $\infty$ has the obvious meaning.

\begin{definition}
A mapping $f: F((G)) \longrightarrow F((G))$ is said to be $v$-compatible if
\[v(m) \leq v(m') \Longleftrightarrow v(f(m)) \leq v(f(m'))\]
for all $m,m' \in F((G))$.
\end{definition}

In \cite{D1} a $v$-compatible mapping is called monotone. For arbitrary $m,m' \in F((G))$ and $a \in F$ 
we define $m+m'$ and $ma$ as usual:
\[m+m': G \longrightarrow F, g \longmapsto m_g+m'_g,\]
\[ma: G \longrightarrow F, g \longmapsto m_ga.\]
These operations are well-defined and they turn $F((G))$ into a right $F$-vector space. The endomorphism 
ring $\End(F((G)))$ of $F((G))$ is a left $F$-vector space with respect to the usual operations, that 
is, $\End(F((G)))$ is a left $F$-algebra.

\begin{proposition}(cf. \cite[Proposition 5.2]{D1})\label{vvertraeglichistinjektiv}
Any $v$-compatible endomorphism of $F((G))$ is injective. The set of all $v$-compatible automorphisms 
of $F((G))$ is a subgroup of $\Aut F((G))$.
\end{proposition}

Later on we will see that there exists a naturally embedding of $F[G,\eta,\alpha]$ into the endomorphism 
ring of $F((G))$. This embedding depends on the basis $X_G$ that has been fixed above and for a different 
basis $X'_G$ of $F[G,\eta,\alpha]$ we will get a different endomorphism ring and therefore another 
embedding of $F[G,\eta,\alpha]$.\smallskip\\
An element $m \in F((G))$ shall also be written as a formal sum $m = \sum x_g m_g$ such that
\[\sum x_gm_g + \sum x_g m'_g = \sum x_g(m_g + m'_g),\ (\sum x_g m_g)a = \sum x_g(m_ga)\]
for all $m,m' \in F((G))$ and $a \in F$. In case of $\supp m = \{g_1,\dots,g_n\}$ we also use the 
notation $m = x_{g_1}m_{g_1} + \dots + x_{g_n}m_{g_n}$. To avoid any misunderstanding we 
occasionally write $F((X_G))$ and $F((X'_G))$ instead of $F((G))$ if $F[G,\eta,\alpha]$ will be 
considered with respect to different bases $X_G$ and $X'_G$ respectively.
\begin{definition}\label{summierbar}
A set $\{m_i \ |\ i \in I\} \subseteq F((G))$ of formal power series is said to be summable if for 
any $g \in G$ there is only a finite number of $i \in I$ such that $m_i(g) \not= 0$ and if the union 
of all corresponding supports $\supp m_i$, $i \in I$ is well-ordered.
\end{definition}
  
If $\{m_i \ |\ i \in I\}$ is summable then $\sum_{i \in I} m_i(g)$ is a finite sum for any $g \in G$ 
such that $m_g := \sum_{i \in I} m_i(g) \in F$ and the support of 
$m: G \longrightarrow F, g \longmapsto m_g$ is well-ordered. In this case we write 
$\sum_{i \in I} m_i = m \in F((G))$ and we obtain a new meaning of the representation of  
formal power series as formal sums. To see this, let $m = \sum x_gm_g$ be a non-zero element of 
$F((G))$. Considering $x_gm_g$ for $g \in \supp m$ as the formal power series with support $\{g\}$ 
which maps $g$ onto $m_g$ we immediately conclude that $\{x_gm_g \ |\ g \in \supp m\}$ is summable 
and
\[\sum_{g \in \supp m} x_gm_g = m.\]

\begin{lemma}\label{summensummen}
Let $I,J$ be sets and for all $i \in I, j \in J$ let $m_{ij} \in F((G))$ such that the following conditions 
are satisfied:
\begin{enumerate}
\item[(1)] $\supp m_{ij}$ and $\supp m_{ij'}$ are disjoint for all $i \in I$ and different $j,j' \in J$,
\item[(2)] $\{m_{ij} \mid j \in J\}$ is summable for all $i \in I$,
\item[(3)] If $m_i = \sum_{j \in J} m_{ij}$ for all $i \in I$ then $\{m_i \mid i \in I\}$ is summable.           
\end{enumerate}

Furthermore, let $K$ be a set and let $M_k \subseteq I \times J$ for all $k \in K$ such that 
$M_k \cap M_{k'} = \emptyset$ for $k \not= k'$ and $I \times J = \bigcup_{k \in K} M_k$. Then, the 
following statements hold true:
\begin{enumerate}
\item[(4)] $\{m_{ij} \mid (i,j) \in M_k\}$ is summable for all $k \in K$,
\item[(5)] If $\om_k = \sum_{(i,j) \in M_k} m_{ij}$ for all $k \in K$ then $\{\om_k \mid k \in K\}$ 
           is summable,
\item[(6)] $\sum_{i \in I} m_i = \sum_{k \in K} \om_k$.
\end{enumerate}
\end{lemma}

{\bf Proof.} We first observe that $M = \bigcup_{i \in I, j \in J} \supp m_{ij}$ is well-ordered. Indeed, 
because of assumption (1) and (2) for all $i \in I$ the support $\supp m_i$ coincides with the union of 
all $\supp m_{ij}$, $j \in J$ and assumption (3) yields the claim.\smallskip\\
Next, we show statement (4). As a subset of $M$, the union of all $\supp m_{ij}$ with $(i,j) \in M_k$ 
is well-ordered. Now, let $g \in G$. We choose $(i_1,j_1), \dots , (i_l, j_l) \in M_k$ where $i_1, \dots ,i_l$ 
are pairwise distinct and $m_{i_1j_1}(g), \dots , m_{i_lj_l}(g)$ are non-zero. Assumption (1) then 
implies $m_{i_1}(g), \dots , m_{i_l}(g) \not= 0$. According to condition (3) there is only a finite 
number of $i \in I$ such that $(i,j) \in M_k$ and $m_{ij}(g) \not= 0$ for some $j \in J$. Let these 
indices be given as $i_1, \dots ,i_l$. Because of assumption (2) for any $i \in \{i_1, \dots ,i_l\}$ 
there is only a finite number of $j \in J$ satisfying $m_{ij}(g) \not= 0$ (condition (1) even ensures 
that there is just one $j$). Hence, the number of all $(i,j) \in M_k$ with $m_{ij}(g) \not= 0$ is 
finite.\smallskip\\
We turn to the proof of statement (5). The same arguments as used above provide that 
$\bigcup_{k \in K} \supp \om_k$ is well-ordered. We assume that there exist $g \in G$ and an infinite 
number of pairwise distinct $k_1,k_2, \dots  \in K$ such that $\om_{k_1}(g), \om_{k_2}(g), \dots $ are 
not-trivial. Then, for any $l=1,2,\dots$ there exists $(i_l,j_l) \in M_{k_l}$ with $m_{i_lj_l}(g) \not= 0$ 
and therefore $m_{i_l}(g) \not= 0$. We apply assumption (3) and conclude that there is an $i \in I$ 
and an infinite number of $j \in J$ satisfying $m_{ij}(g) \not= 0$. This contradicts condition (2).
\smallskip\\
In order to verify statement (6), let $g \in G$ with $\{i \in I \mid m_i(g) \not= 0\} = 
\{i_1, \dots , i_k\}$. By assumption (1) for any $l = 1, \dots ,k$ there exists exactly one $j_l \in J$ 
with $m_{i_l j_l}(g) \not= 0$, that is, $m_{i_l}(g) =  m_{i_l j_l}(g)$. Thus, 
$m_{i_1 j_1}(g) + \dots + m_{i_k j_k}(g)$ is the image of $g$ with respect to $\sum_{i \in I} m_i$. 
On the other hand, for any $k \in K$ with $\om_k(g) \not= 0$ there exists $(i,j) \in M_k$ satisfying 
$m_{ij}(g) \not= 0$. Condition (1) implies $m_i(g) \not= 0$  which means $i \in \{i_1, \dots , i_k\}$ 
and $j \in \{j_1, \dots , j_k\}$. Hence, $\om_k(g)$ is the sum of all $m_{i_l j_l}(g)$ such that 
$(i_l, j_l) \in M_k$ and therefore $m_{i_1 j_1}(g) + \dots + m_{i_k j_k}(g)$ is also 
the image of $g$ with respect to $\sum_{k \in K} \om_k$.\\
\qed\\

\begin{definition}
An endomorphism $f$ of $F((G))$ is called continuous if the following condition is satisfied: 
If $\{m_i \ |\ i \in I\} \subseteq F((G))$ is summable then $\{f(m_i) \ |\ i \in I\}$ is summable and 
$f(\sum_{i \in I} m_i) = \sum_{i \in I} f(m_i)$. 
\end{definition}

Continuous endomorphisms are called $\sigma$-linear in \cite{D1}.

\begin{proposition} (cf. \cite[Proposition 5.1]{D1})\label{speziellesstetig}
An endomorphism $f$ of $F((G))$ is continuous if for any $m = \sum x_gm_g \in F((G))$ 
the set $\{ f(x_gm_g) \ |\ g \in \supp m\}$ is summable and $f(m) = \sum_{g \in \supp m} f(x_gm_g)$. 
The set of all continuous endomorphisms of $F((G))$ is a subalgebra of $\End (F((G)))$.
\end{proposition}

\begin{theorem} (cf. \cite[Lemma 3]{D3} and \cite{S})\label{finvstetig}
If $f$ is a continuous and $v$-compatible automorphism of $F((G))$ then $f\inv$ is also continuous.
\end{theorem}

\begin{definition}
A mapping $f: F((G)) \longrightarrow F((G))$ is called $v$-compatible on $G$ if
\[g \leq g' \Longleftrightarrow v(f(x_g)) \leq v(f(x_{g'}))\]
for all $g,g' \in G$. And $f$ is called surjective on $G$ if for any $g \in G$ there exists $g' \in G$ 
such that $g = v(f(x_{g'}))$.
\end{definition}

In \cite{D1} the properties {\it $v$-compatible on $G$} and {\it surjective on $G$} are called 
{\it locally monotone} and {\it locally surjective} respectively.

\begin{theorem}(cf. \cite[Theorem 5.1]{D1})\label{lokalglobal}
A continuous endomorphism $f$ of $F((G))$ is a $v$-compatible automorphism if and only if $f$ is surjective 
on $G$ and $v$-compatible on $G$.
\end{theorem}

Next, we explain how $R=F[G,\eta,\alpha]$ can be considered as a subring of the endomorphism ring 
of $F((G))$ where the assumptions are as stated above (further details can be found in \cite{G,GS,J}). 
At first, we observe that $F[G,\eta,\alpha]$ can be considered as a subgroup of $F((G))$ with respect 
to the addition since any $m = a_{g_1}x_{g_1} + \dots + a_{g_n}x_{g_n} \in F[G,\eta,\alpha]$ can be written 
as $m = x_{g_1}m_{g_1} + \dots + x_{g_n}m_{g_n}$ with $m_{g_1}, \dots ,m_{g_n} \in F$ which causes an 
isomorphism between $F[G,\eta,\alpha]$ and the additive group of all formal power series with finite 
support. Occasionally, it is convenient to interpret power series with finite supports as elements 
from $F[G,\eta,\alpha]$ in order to apply the multiplication given in $F[G,\eta,\alpha]$. We will 
use this frequently without any further explanation.\smallskip\\
Let us return to the problem of embedding $R=F[G,\eta,\alpha]$ into $\End(F((G)))$. If $m=\sum x_gm_g$ 
is an arbitrary element of $F((G))$ then $\{hg \ |\ g \in \supp m\}$ is well-ordered for any $h \in G$ 
since $\leq$ is a left-order of $G$. This ensures that
\[\iota_{a_h x_h}: F((G)) \longrightarrow F((G)), \sum x_gm_g \longmapsto \sum x_{hg}(\alpha\inv_{hg}(a_h\eta(h,g))m_g)\]
is a well-defined mapping. According to the remark above, it should always be noted that  
\[x_{hg}(\alpha\inv_{hg}(a_h\eta(h,g))m_g) = a_h x_h \gap x_gm_g\]
for arbitrary $g,h \in G$ and $m_g, a_h \in F$ which yields 
\[\iota_{a_h x_h} \circ \iota_{a'_{h'} x_{h'}} = \iota_{a_h x_h a'_{h'} x_{h'}}\]
for all $h,h' \in G$ and $a_h,a'_{h'} \in F$. If $a_h \not=0$ then $\iota_{a_hx_h}$ is a 
continuous and $v$-compatible automorphism of $F((G))$ with $\iota_{a_h x_h}(x_g F) \subseteq x_{hg}F$ 
for all $g \in G$. Following \cite{D1,D2}, we call endomorphisms with this property monomial. 
For pairwise distinct $h_1,\dots,h_n \in G$ and arbitrary $a_{h_1}, \dots ,a_{h_n} \in F$ we now 
define
\[\iota_{a_{h_1}x_{h_1} + \dots + a_{h_n}x_{h_n}}: F((G)) \longrightarrow F((G)), 
         m \longmapsto \iota_{a_{h_1}x_{h_1}}(m) + \dots + \iota_{a_{h_n}x_{h_n}}(m).\]
Clearly, $\iota_{a_{h_1}x_{h_1} + \dots + a_{h_n}x_{h_n}}$ is an endomorphism of $F((G))$ which in 
addition is continuous because of Proposition \ref{speziellesstetig}. A straightforward calculation 
shows that
\[\iota: F[G,\eta,\alpha] \longrightarrow \End(F((G))),
         a_{h_1}x_{h_1} + \dots + a_{h_n}x_{h_n} \longmapsto \iota_{a_{h_1}x_{h_1} + \dots + a_{h_n}x_{h_n}}\]
is an injective ring homomorphism. Hence, with respect to the composition $\iota(\G)$ is a group of 
monomial, $v$-compatible, and continuous automorphisms of $F((G))$.\smallskip\\
For the sake of simplicity, we shall not distinguish between the elements of $F[G,\eta,\alpha]$ and 
their images under $\iota$. It will be clear from the context whether 
$a_{h_1}x_{h_1} + \dots + a_{h_n}x_{h_n}$ is an element from $F[G,\eta,\alpha]$ or $F((G))$ or 
whether it is the endomorphism of $F((G))$ which maps $m = \sum x_gm_g$ onto 
$a_{h_1}x_{h_1}m + \dots + a_{h_n}x_{h_n}m$ where $a_{h_i}x_{h_i}m = \sum a_{h_i}x_{h_i}x_gm_g$ for 
$i=1,\dots,n$. Therefore, we also consider $\G$ as a group of monomial, $v$-compatible, and continuous 
automorphisms of $F((G))$.\\

\begin{proposition} (cf. \cite{D0,D1,D2}) \label{ohneadditiverelationen}
Let the notation be as above. Then, any non-zero element from $F[G,\eta,\alpha]$ - regarded as an 
endomorphism of $F((G))$ - is a continuous and $v$-compatible automorphism.
\end{proposition}

{\bf Proof.} We apply Proposition 8.1 of \cite{D1}. Let $\G_h = \{a_gx_g \in \G \ |\ v(a_gx_gx_h) = h\}$ 
where $h \in G$. Clearly, $v(a_gx_gx_h)=gh$ implies $\G_h = \{ax_e \ |\ a \in F^\times\}$ and we obtain
for $a_1x_e, \dots ,a_nx_e \in \G_h$ with
\[\Null = (a_1x_e + \dots + a_nx_e)(x_h) = a_1x_ex_h + \dots a_nx_ex_h = (a_1 + \dots + a_n)x_h\]
that $a_1 + \dots + a_n = 0$ and therefore $a_1x_e + \dots + a_nx_e = \Null$. This proves the claim 
since $F[G,\eta,\alpha]$ is the subring of $\End(F((G)))$ generated by $\G$.\\
\qed\\

\begin{definition}
Let $F[G,\eta,\alpha]$ be a crossed product of $G$ over the skew field $F$ with respect to the basis $X_G$ 
and let $\leq$ be a left-order of $G$. The rational closure $F(G,\eta,\alpha)$ of $F[G,\eta,\alpha]$ 
in $\End(F((G)))$ (the rational closure of $\G$ in $\End(F((G)))$) is called the Dubrovin-quotient 
ring of $F[G,\eta,\alpha]$ (with respect to $\leq$ and $X_G$). If $F(G,\eta,\alpha)$ is a division 
ring then $F(G,\eta,\alpha)$ is called the Dubrovin-division ring of fractions of $F[G,\eta,\alpha]$.
\end{definition}

\begin{proposition}\label{fortsetzungenvonautomorphismen}
Let $R=F[G,\eta,\alpha]$ be a crossed product of $G$ over the skew field $F$ with respect to the basis $X_G$ 
and let $\leq$ be a left-order of $G$. Further, let $\varphi: R \longrightarrow R$ be a ring automorphism.  
If $D$ is the Dubrovin-quotient ring of $F[G,\eta,\alpha]$ with respect to $\leq$ and $X_G$ and if $D'$ 
is the Dubrovin-quotient ring of $F[G,\eta,\alpha]$ with respect to $\leq$ and $X'_G = \{x'_g \ |\ g \in G\}$ 
where $x'_g = \varphi(x_g)$ for all $g \in G$ then there exists a unique ring isomorphism 
$\psi: D \longrightarrow D'$ which extends $\varphi$.
\end{proposition}

{\bf Proof.} The uniqueness of $\psi$ can easily be proved by transfinite induction on the complexity 
$cp(d)$, $d \in D$ which describes how $d$ is built in $\End(F((G)))$ from elements of $\G$. Since 
$\varphi(F)=F$ by Proposition \ref{Einheitengruppe} we obtain that
\[\rho: F((X_G)) \longrightarrow F((X'_G)), 
  \sum x_gm_g \longmapsto \sum x'_g \varphi(m_g) = \sum \varphi(x_g) \varphi(m_g)\]
is a group isomorphism between the addivite groups $F((X_G))$ and $F((X'_G))$ such that 
$\rho(ma) = \rho(m)\varphi(a)$ for all $m \in F((X_G))$ and $a \in F$. A straightforward calculation 
shows that
\[\tau: \End (F((X_G))) \longrightarrow \End (F((X'_G))), f \longmapsto \rho \gap f \gap \rho\inv \]
is a ring isomorphism such that $\tau(a_hx_h) = \varphi(a_h x_h)$ for all $h \in G$ and $a_h \in F$. 
In addition, $\tau$ and $\varphi$ coincide on $F[G,\eta,\alpha]$ and therefore $\tau$ maps the 
Dubrovin-quotient ring $D$ of $F[G,\eta,\alpha]$ with respect to $\leq$ and $X_G$ onto the Dubrovin-quotient 
ring $D'$ of $F[G,\eta,\alpha]$ with respect to $\leq$ and $X'_G$. The restriction $\psi$ of $\tau$ 
to $D$ is the desired isomorphism.\\
\qed\\

The following lemma occurs in \cite[Proposition 7.1]{D1,D2} in a slightly modified version. For the 
sake of completeness we add Dubrovin's proof.

\begin{lemma}\label{duaufxeanwenden}
Let $R=F[G,\eta,\alpha]$ be a crossed product of the left-ordered group $G$ over the skew field $F$ 
and let $U$ be a subgroup of $G$. If $D_U$ denotes the rational closure of $F[U,\eta,\alpha]$ in 
$\End(F((G)))$ then 
\[\supp f(m) \subseteq \{gh \ |\ g \in U, h \in \supp m\}\] 
for all $f \in D_U$ and $m \in F((G))$.
\end{lemma}

{\bf Proof.} Let the notation be as above, especially let $\End(F((G)))$ be defined with respect to 
the fixed basis $X_G$ of $F[G, \eta, \alpha]$. We prove the lemma by transfinite induction on the 
complexity $cp(f)$, $f \in D_U$ which describes how $f$ is built in $\End(F((G)))$ from elements of 
$\U = F^\times X_U$. The statement is obvious in case of $f=\Null$ or $f \in \U$. Thus, let $cp(f)>1$. 
If $f$ is additively decomposable then there exist $f_1,\dots,f_n \in D_U$ such that $f_1,\dots,f_n \lhd f$ 
and $f = f_1 + \dots + f_n$. Now, the claim follows immediately by 
$\supp f(m) \subseteq \supp f_1(m) \cup \dots \cup \supp f_n(m)$ and the induction hypothesis. Next, 
let $f$ be additively indecomposable. Then, $f$ is a unit in $\End (F((G)))$ by \cite[Theorem 4.9]{DGH} 
and therefore in $D_U$ where $f\inv \lhd f$. There exist uniquely determined $m_1, m_2 \in F((G))$ with 
$f(m) = m_1 + m_2$ such that
\[\supp m_1 \subseteq \{gh\ |\ g \in U, h \in \supp m\} \ \mbox{and}\  
 \supp m_2 \cap \{gh\ |\ g \in U, h \in \supp m\} = \emptyset\]
where $m_2 = \Null$ must be shown. Clearly, $f\inv(m_2) = m - f\inv(m_1)$ and the induction hypothesis 
imply $\supp f\inv (m_2) \subseteq \{gh\ |\ g \in U, h \in \supp m\}$. Now, let us assume that 
$m_2 \not= \Null$. Then, $\supp f\inv (m_2) \not= \emptyset$ and therefore
\[\{gh\ |\ g \in U, h \in \supp m_2\} \cap \{gh\ |\ g \in U, h \in \supp m\} \not= \emptyset.\]
Hence, there are $g,g' \in U$, $h \in \supp m_2$, and $h' \in \supp m$ such that $gh=g'h'$. This yields 
the contradiction
\[\supp m_2 \cap \{gh\ |\ g \in U, h \in \supp m\} \not= \emptyset.\]
\qed\\

The next corollary is a direct consequence of Dubrovin's lemma and appears in \cite{J}.

\begin{corollary}\label{dubrovinisthughes}
Let $R=F[G,\eta,\alpha]$ be a crossed product of a left-ordered group $G$ over a skew field $F$ 
with respect to the basis $X_G$ and let $F(G,\eta,\alpha)$ be the Dubrovin-quotient ring of 
$F[G,\eta,\alpha]$ with respect to $\leq$ and $X_G$. If $F(G,\eta,\alpha)$ is a division ring then 
$F(G,\eta,\alpha)$ is strongly Huges-free.
\end{corollary}

{\bf Proof.} Let $H$ be an arbitrary subgroup of $G$ with normal subgroup $N$ and let $h_1,\dots,h_n \in H$, 
$n \in \N$ are given such that $h_1N, \dots, h_nN$ are pairwise different. It has to be shown that the 
mappings $x_{h_1}, \dots , x_{h_n}$ are linearly independent over the rational closure $D_N$ of 
$F[N,\eta,\alpha]$ in $F(G,\eta,\alpha)$. Thus, let $d_1x_{h_1} + \dots + d_nx_{h_n}$ be the zero-map 
for arbitrarily given $d_1, \dots ,d_n \in D_N$. Then $d_1x_{h_1}(x_e) + \dots + d_nx_{h_n}(x_e) = \Null$ 
which means $d_1(x_{h_1}) + \dots + d_n(x_{h_n}) = \Null$. Since $Nh_1, \dots , Nh_n$ are pairwise disjoint 
subsets of $G$ we conclude that $d_i(x_{h_i})$, $i=1,\dots,n$ have disjoint supports by Lemma 
\ref{duaufxeanwenden}. Thus, $d_i(x_{h_i}) = \Null$ and therefore $d_i = \Null$ for all $i=1,\dots,n$.\\
\qed\\

\section{Free Division Rings of Fractions and Power Series}\label{Quotientenschiefkoerper und Potenzreihen}

\vspace*{0.4cm}

We shall use the notation as introduced in the sections before.

\begin{theorem} \label{Fortsetzungvonphi}
Let $R=F[G,\eta,\alpha]$ be a crossed product of a locally indicable group $G$ over a skew field 
$F$ with respect to the basis $X_G$ and let $\leq$ be a Conradian left-order of $G$. Furthermore, 
let $D$ be a division ring of fractions of $R$ which is free with respect to $\leq$ and let $\R$ be 
the Dubrovin-quotient ring of $F[G,\eta,\alpha]$ with respect to $\leq$ and $X_G$. Then, there exists 
a ring isomorphism $\psi: D \longrightarrow \R$ which is the identity on $F[G,\eta,\alpha]$. Especially, 
$\R$ is a division ring and $D$ is strongly Hughes-free.
\end{theorem}

\begin{corollary} \label{alleisomorph}
Let $R=F[G,\eta,\alpha]$ be a crossed product of a locally indicable group $G$ over a skew field 
$F$ with respect to the basis $X_G$ and let $D$,$D'$ be two division rings of fractions of $R$ which 
are free with respect to two Conradian left-orders $\leq$,$\leq'$ respectively. Any ring automorphism 
$\varphi$ of $R$ can be extended to a ring isomorphism $\psi: D \longrightarrow D'$.
\end{corollary}

{\bf Proof}. By Theorem \ref{Fortsetzungvonphi} we may assume that $D=\R$ where $\R$ is the 
Dubrovin-division ring of fractions of $F[G,\eta,\alpha]$ with respect to $\leq$ and $X_G$. Furthermore, 
by Theorem \ref{Fortsetzungvonphi}, Propositon \ref{hughes-free ist free}, and 
Proposition \ref{freifuervarphi(R)} we may also assume that $D' = \R'$ where $\R'$ is the Dubrovin-division 
ring of fractions of $F[G,\eta,\alpha]$ with respect to $\leq$ and $X'_G = \varphi(X_G)$. It remains 
to apply Proposition \ref{fortsetzungenvonautomorphismen}.\\
\qed\\

{\bf Remark.} Theorem \ref{Fortsetzungvonphi} is proved in \cite{J} under the special assumption that 
$\leq$ has maximal rank. Furthermore, in \cite{J} it is also shown that Hughes-freeness and strongly 
Hughes-freeness coincide which solves P.A. Linnell's Problem 4.8 from (cf. \cite{L}). Here it is a
consequence of Theorem \ref{lokalindizierbaristconrad}, Proposition \ref{existenzmaximalerrang}, 
Proposition \ref{hughes-free ist free}, and Theorem \ref{Fortsetzungvonphi}. Our results also provide 
the following

\begin{corollary}
Let $R=F[G,\eta,\alpha]$ be a crossed product of a locally indicable group $G$ over a skew field 
$F$ and let $D$ be a division ring of fractions of $R$. Then $D$ is (strongly) Hughes-free 
if and only if $D$ is free with respect to at least one Conradian left-order of $G$. In this case, 
$D$ is free with respect to any Conradian left-order of $G$.
\end{corollary}
By means of \cite{JL} we also obtain

\begin{corollary}
Let $G$ be a locally indicable group and let $F$ be a commutative field of characteristic zero. Then 
for any Conradian left-order the Dubrovin-quotient ring of the group ring $F[G]$ is a division ring.
\end{corollary}

As we will explain at the end of this section, Theorem \ref{Fortsetzungvonphi} is an immediate 
consequence of Lemma \ref{aufbauvonpsi} below. The notation we use is the same as in 
Theorem \ref{Fortsetzungvonphi} and Section \ref{Komplexitaet} where Definition \ref{lambdausw} is 
of particular importance. The complexity $cp(d)$ for some $d \in D$ always refers to the manner how 
$d$ is built from elements of $\G = F^\times G$.

\begin{lemma}\label{aufbauvonpsiallgemeiner}
Let $\lambda \in \Lambda$ and let $\psi_\lambda: D_{\leq \lambda} \longrightarrow \R$ be a mapping 
such that the following conditions are fulfilled:
\begin{enumerate}
\item[(1)] $\psi_\lambda(x)=x$ for all $x \in F[G,\eta,\alpha]$.
\item[(2)] $\psi_\lambda(x)$ is a continuous endomorphism of $F((G))$ for all 
           $x \in D_{\leq \lambda}$.
\item[(3)] If $x \in D_{\leq \lambda}$ is reduced with $cp(x) > 1$ and if 
           $x = \sum _{g \in \fC} d_g \gap x_g$ is the pure left representation of $x$ then 
           $(\psi_\lambda(x))(x_e)=\sum _{g \in \fC} (\psi_\lambda(d_g) \gap x_g)(x_e) \in F((C))$. 
\item[(4)] $\psi_\lambda(x_\tg \gap x) = x_\tg \gap \psi_\lambda(x)$ and 
           $\psi_\lambda(x \gap x_\tg) = \psi_\lambda(x) \gap x_\tg$ for all $x \in D_{\leq \lambda}$ 
           and all $\tg \in G$.
\end{enumerate}

Then the following hold true:

\begin{enumerate}
\item[(5)] If $x \in D_{\leq \lambda}$, $cp(x) > 1$ is arbitrary and if 
           $x_h^{-1} \gap x = \sum _{g \in \fC} d_g \gap x_g$ is an arbitrary left representation of $x$ then
           $(\psi_\lambda(x))(x_e) = x_h \gap \sum _{g \in \fC} (\psi_\lambda(d_g) \gap x_g)(x_e) \in x_hF((C))$.
\item[(6)] If $H$ is a convex subgroup of $G$ and $D_H$ the rational closure of $F[H,\eta,\alpha]$ 
           in $D$ then $(\psi_\lambda(x))(x_e) \in F((H))$ for any $x \in D_{\leq\lambda}$, $cp(x)>1$ which is 
           in $D_H$. Especially, $(\psi_\lambda(d_g))(x_e) \in F((C'))$ for all $d_g$ occurring 
           as a coefficient of a left representation of $x$ as given in (5).
\item[(7)] If $x \in D_{\leq \lambda}$, $cp(x) > 1$ is arbitrary and if 
           $x_h^{-1} \gap x = \sum _{g \in \fC} d_g \gap x_g$ is an arbitrary left representation of $x$ then
           $(\psi_\lambda(x))(x_\og) = x_h \gap \sum _{g \in \fC} (\psi_\lambda(d_g) x_g)(x_\og) \in x_hF((C))$  
           for all $\og \in C$.
\item[(8)] If $x \in D_{\leq \lambda}$, $cp(x) > 1$ is arbitrary and if 
           $x_h^{-1} \gap x = \sum _{g \in \fC} d_g \gap x_g$ is an arbitrary left representation of $x$ then
           $(\psi_\lambda(x))(m) = x_h \gap \sum _{g \in \fC} (\psi_\lambda(d_g) x_g)(m) \in x_hF((C))$  
           for all $m \in F((C))$.
\end{enumerate}
\end{lemma}
{\bf Remarks.}
\begin{enumerate}
\item[(1)] Assumption (3) includes that $\{(\psi_\lambda(d_g) \gap x_g)(x_e) \mid g \in \fC\}$ is 
           summable.
\item[(2)] The appearance of an infinite sum in the statements above always include that this sum 
           actually exists.
\item[(3)] Proposition \ref{Dleqlambda} and $d_g \lhd x$ imply $d_g \in D_{\leq \lambda}$ for any $d_g$ 
           which occurs as a coefficient of a left representation of $x$ as given above.         
\end{enumerate}

{\bf Proof.} In order to prove statement (5) we first observe that $x\inv_h x$ is reduced with pure 
left representation $x_h^{-1} \gap x = \sum _{g \in \fC} d_g \gap x_g$. Clearly, 
$x\inv_h x \in D_{\leq \lambda}$ and therefore 
$(\psi_\lambda(x\inv_h x))(x_e) = \sum _{g \in \fC} (\psi_\lambda(d_g) \gap x_g)(x_e) \in F((C))$ 
by assumption (3). We now apply assumption (4) and conclude 
$(\psi_\lambda(x))(x_e) = (\psi_\lambda(x_h(x\inv_h x)))(x_e) = x_h (\psi_\lambda(x\inv_h x))(x_e)$ 
which completes the proof of statement (5). Furthermore, statement (6) follows from statement (5) and 
remark (2) after Theorem \ref{EindeutigkeitderDarstellung}. We show statement (7). For any $g \in \fC$ 
there are $g' \in \fC$ and $u_g \in F^\times X_{C'}$ with $x_g x_\og = u_g x_{g'}$. Since 
$g \og C' = \og g C' = g'C'$ for all $g \in \fC$ the set 
$\{g' \mid g \in \fC \mbox{\ and\ } d_g \not= \Null \}$ is well-ordered with respect to $\leq$. This 
means that $x\inv_h(x x_\og) = \sum_{g \in \fC} (d_g u_g) x_{g'}$ is a left representation of 
$x x_\og \in D_{\leq \lambda}$ and statement (5) yields 
$(\psi_\lambda(x x_\og))(x_e) = x_h \sum_{g \in \fC} (\psi_\lambda(d_g u_g) x_{g'})(x_e) \in x_h F((C))$. 
Assumption (4) provides $(\psi_\lambda(d_g u_g)x_{g'})(x_e) = (\psi_\lambda(d_g u_g x_{g'}))(x_e) 
= (\psi_\lambda(d_g) x_g)(x_{\og})$ and $(\psi_\lambda(x x_\og))(x_e)=(\psi_\lambda(x))(x_\og)$ 
which proves statement (7). In order to derive statement (8) let  
$m = \sum_{\og \in C} x_\og \gap m_\og$ be an arbitrary element of $F((C))$. Assumption (2) 
ensures that $\psi_\lambda(x)$ is continuous and therefore 
$(\psi_\lambda(x))(m) = \sum_{\og \in C} (\psi_\lambda(x))(x_\og \gap m_\og)$. Thus, 
\[(\psi_\lambda(x))(m) = \sum_{\og \in C} x_h\sum _{g \in \fC} (\psi_\lambda(d_g) \gap x_g)(x_\og \gap m_\og)\]
because of statement (7). By means of $d_g \in D_{\leq\lambda}$ and assumption (4) as well as 
$d_g \in D^-_C$ and statement (6) we conclude
\[(\psi_\lambda(d_g) \gap x_g)(x_\og \gap m_\og) = 
            [x_g x_\og (\psi_\lambda( x_\og^{-1} x_g^{-1}\gap d_g\gap x_g x_\og))(x_e)] m_\og \in x_g x_\og\gap F((C')).\]
This shows that for a fixed $\og \in C$ the sets $\supp (\psi_\lambda(d_g) \gap x_g)(x_\og\gap m_\og)$, 
$g \in \fC$ are pairwise disjoint. We apply Lemma \ref{summensummen} and obtain

\begin{eqnarray*}          
          (\psi_\lambda(x))(m) &=& \sum_{\og \in C} x_h\sum _{g \in \fC} (\psi_\lambda(d_g) \gap x_g)(x_\og \gap m_\og) 
                                =  x_h\sum _{g \in \fC} \sum_{\og \in C} (\psi_\lambda(d_g) \gap x_g)(x_\og \gap m_\og)\\
                               &=& x_h\sum _{g \in \fC} (\psi_\lambda(d_g) \gap x_g)(\sum_{\og \in C} x_\og \gap m_\og)
                                =  x_h\sum _{g \in \fC} (\psi_\lambda(d_g) \gap x_g)(m).
\end{eqnarray*}
As we have seen above, $(\psi_\lambda(d_g) \gap x_g)(x_\og \gap m_\og) \in x_g x_\og\gap F((C')) \subseteq F((C))$ 
for all $g \in \fC$, $\og \in C$ and we conclude $(\psi_\lambda(x))(m) \in x_hF((C))$ which completes the proof.\\
\qed\\

Now, we are in the position to prove

\begin{lemma}\label{aufbauvonpsi}
For all $\lambda \in \Lambda$ there exists a mapping  
$\psi_\lambda: D_{\leq \lambda} \longrightarrow \R$ such that the following hold true:
\begin{enumerate}
\item[(1)]  If $\lambda' \in \Lambda$ with $\lambda' < \lambda$ then $\psi_\lambda(x) = \psi_{\lambda'}(x)$ 
            for all $x \in D_{\leq\lambda'}$.
\item[(2)]  $\psi_\lambda(x)=x$ for all $x \in F[G,\eta,\alpha]$.
\item[(3)]  $\psi_\lambda$ is additive on $D_{< \lambda}$.
\item[(4)]  If $a \in D_{\leq \lambda}$ is a proper atom and $\lambda' = cp(a)$ then the restriction 
            of $\psi_\lambda$ to $D_{< \lambda'}$ is an injective ring homomorphism.
\item[(5)]  $\psi_\lambda(x)$ is a continuous and $v$-compatible automorphism of $F((G))$ for all 
            non-zero $x \in D_{\leq \lambda}$.
\item[(6)]  If $x \in D_{\leq \lambda}$ is reduced with $cp(x) > 1$ and if 
            $x = \sum _{g \in \fC} d_g \gap x_g$ is the pure left representation of $x$ then 
            $(\psi_\lambda(x))(x_e)=\sum _{g \in \fC} (\psi_\lambda(d_g) \gap x_g)(x_e) \in F((C))$. 
\item[(7)]  $\psi_\lambda(x_\tg \gap x) = x_\tg \gap \psi_\lambda(x)$ and 
            $\psi_\lambda(x \gap x_\tg) = \psi_\lambda(x) \gap x_\tg$ for all $x \in D_{\leq \lambda}$ 
            and $\tg \in G$.
\end{enumerate}
\end{lemma}

{\bf Proof.} We prove Lemma \ref{aufbauvonpsi} by transfinite induction on $\lambda \in \Lambda$ and 
consider the case $\lambda = 1$ first. Then, $D_{\leq \lambda} = F[G,\eta,\alpha]$ and we choose $\psi_1$ 
as the identity map such that statements (1),(2),(3),(4), and (7) are obviously true. Clearly, statement (5) 
follows by Proposition \ref{ohneadditiverelationen}. To prove statement (6) we first choose a left 
representation of $x$ as given in the proof of Theorem \ref{DarstellungeninfreienKoerpern} which is a 
finite sum $x_h^{-1}x = \sum_{g \in \fC} d_g \gap x_g$ with $h \in C$ since $x$ is reduced and 
$d_g \in F[C',\eta,\alpha]$ for all $g \in \fC$. We now transform this left representation into the 
pure left representation of $x$ as explained before. Thus $h=e$ can be assumed, that is, 
$\psi_1(x)(x_e) = x = \sum_{g \in \fC} d_g \gap x_g 
                  = \sum_{g \in \fC} (\psi_1(d_g) \gap x_g)(x_e) \in F((C))$.\smallskip\\
Now, let $x \in D_{\leq\lambda}$ where $\lambda > 1$. If there exists $\lambda' \in \Lambda$ with 
$\lambda' < \lambda$ such that $x \in D_{\leq \lambda'}$ then we define 
$\psi_\lambda (x)=\psi_{\lambda'}(x)$. Clearly, this definition does not depend on $\lambda'$, 
$\lambda' < \lambda$ and statements (1) and (2) are obviously true. Even though we have defined 
$\psi_\lambda(x)$ just in this special case we are already in the position to prove statement (3). 
Thus, let $x,y \in D_{< \lambda}$. Because of Proposition \ref{Dleqlambda} there exists 
$\lambda' \in \Lambda$ with $\lambda' < \lambda$ such that $x,y,x+y \in D_{\leq \lambda'}$. 
Therefore, $\psi_\lambda(x), \psi_\lambda(y)$, and $\psi_\lambda(x+y)$ are defined and we can apply the 
induction hypothesis to $\lambda'$. Then, conditions (1)-(4) of Lemma \ref{aufbauvonpsiallgemeiner} 
are fulfilled where $\lambda$ is replaced by $\lambda'$ and statements (5)-(8) of 
Lemma \ref{aufbauvonpsiallgemeiner} can be used for $\lambda'$. With this in mind we now prove:\smallskip\\
{\bf A:} $\psi_\lambda(x+y) = \psi_\lambda(x)+\psi_\lambda(y)$ for all $x,y \in D_{< \lambda}$.\smallskip\\
We assume that there exist $x,y \in D_{< \lambda}$ such that 
$\psi_\lambda(x+y) \not= \psi_\lambda(x)+\psi_\lambda(y)$. Let $m^+$ be the larger of the two ordinals 
$cp(x), cp(y)$ and let $m^-$ be the smaller. We shall also assume that among all these $x,y$ the pair 
$(m^+,m^-)$ is minimal with respect to the lexicographical order. This means, that 
$\psi_\lambda(x'+y') = \psi_\lambda(x')+\psi_\lambda(y')$ for all $x',y' \in D$ with 
$x' \unlhd x, y' \unlhd y$ whereby $x' \lhd x$ or $y' \lhd y$. Since $\psi_\lambda(x)$, $\psi_\lambda(y)$, 
and $\psi_\lambda(x+y)$ are continuous we conclude 
$(\psi_\lambda(x))(x_g) + (\psi_\lambda(y))(x_g) \not= (\psi_\lambda(x+y))(x_g)$ for at least one $g \in G$. 
Because of \cite[Proposition 4.8]{DGH} and 
$(\psi_\lambda(x \gap x_g))(x_e) + (\psi_\lambda(y \gap x_g))(x_e) \not= (\psi_\lambda(x \gap x_g + y \gap x_g))(x_e)$ 
we assume $g = e$ and discuss the three cases $\rm A1$, $\rm A2.1$, and $\rm A2.2$ as 
presented in the remarks about left representations of sums and products after Definition \ref{Definitionreduziert}. 
Using the corresponding notation and assumptions we  make the following observations: \smallskip\\
${\bf A1:}$ Then 
$(\psi_\lambda(x_h^{-1} \gap (x+y)))(x_e) = (\psi_\lambda(x_h^{-1} \gap x))(x_e) 
                                          + (\psi_\lambda(x_h^{-1} \gap y \gap x_g^{-1})x_g)(x_e)$ 
and therefore $(\psi_\lambda(x+y))(x_e) = (\psi_\lambda(x))(x_e) + (\psi_\lambda(y))(x_e)$.\smallskip\\
${\bf A2.1:}$ Then 
$(\psi_\lambda(x_h^{-1}(x+y)))(x_e) = \sum_{g \in \fC} (\psi_\lambda(d_g + d_g') \gap x_g)(x_e)$. 
Because of $d_g \lhd x$ and $d_g' \unlhd y$ we now conclude  
$\psi_\lambda(d_g + d_g') \gap x_g = \psi_\lambda(d_g) \gap x_g + \psi_\lambda(d_g') \gap x_g$. 
According to $(\psi_\lambda(x_h^{-1} \gap x))(x_e) = \sum_{g \in \fC} (\psi_\lambda(d_g) \gap x_g)(x_e)$ 
and $(\psi_\lambda(x_h^{-1} \gap y))(x_e) = \sum_{g \in \fC} (\psi_\lambda(d_g') \gap x_g)(x_e)$ 
this means 
$(\psi_\lambda(x_h^{-1}(x+y)))(x_e) = (\psi_\lambda(x_h^{-1} \gap x))(x_e) 
                                    + (\psi_\lambda(x_h^{-1} \gap y))(x_e)$ and therefore
$(\psi_\lambda(x+y))(x_e) = (\psi_\lambda(x))(x_e) + (\psi_\lambda(y))(x_e)$.\smallskip\\
$\bf A2.2:$ Then 
$(\psi_\lambda(x_h^{-1}(x+y)))(x_e) = (\psi_\lambda(d_g + d_g') \gap x_g)(x_e)$. For all $g' \in \fC$, 
$g' \not= g$ we have $d_{g'} + d_{g'}' = 0$, that is 
$\psi_\lambda(d_{g'}) \gap x_{g'} + \psi_\lambda(d_{g'}') \gap x_{g'} = 0$ and 
we can proceed as in case $\rm A2.1$.\smallskip\\
In any case we deduce a contradiction and we turn to statement (4). Clearly, Proposition \ref{Dleqlambda} 
yields $\lambda' \leq \lambda$ which implies 
$D_{< \lambda'} \subseteq D_{< \lambda} \subseteq D_{\leq \lambda}$ where $D_{< \lambda'}$ 
is a subring of $D_{\leq \lambda}$ by \cite[Proposition 4.1]{DGH}. Statement (3) ensures that $\psi_\lambda$ 
is additive on $D_{< \lambda'}$. In order to verify the injectivity it is therefore enough to prove that 
$\psi_\lambda(x)$ is not the zero map for any non-zero $x \in D_{< \lambda'}$. But this follows from 
statement (5) by the induction hypothesis since $x \in D_{\leq \lambda^\ast}$ for some 
$\lambda^\ast < \lambda' \leq \lambda$ by Proposition \ref{Dleqlambda}. It remains to show the 
multiplicativity. The main argument we shall apply is of crucial importance and later on it will 
be used again several times. The idea is based on Proposition \ref{Dleqlambda} which states that 
for a finite number of elements $a_1, \dots ,a_n \in D_{< \lambda'}$ any finite number of sums of 
products of them belongs to a common $D_{\leq \lambda^\ast}$ where $\lambda^\ast < \lambda' \leq \lambda$. 
Hence, for all these terms the corresponding images under $\psi_\lambda$ are already defined and the 
induction hypothesis can be applied. Taking into account all this we are now going to prove:\smallskip\\
{\bf M:} $\psi_\lambda(xy) = \psi_\lambda(x)\psi_\lambda(y)$ for all $x,y \in D_{< \lambda}$.\smallskip\\
We assume that there exist $x,y \in D_{< \lambda}$ such that 
$\psi_\lambda(xy) \not= \psi_\lambda(x)\psi_\lambda(y)$. Then $cp(x), cp(y) > 1$. Similar to the proof 
of statement {\bf A} we can also assume that $x,y$ are minimal in the following sense: If $x',y' \in D$ 
such that $x' \unlhd x, y' \unlhd y$ whereby $x' \lhd x$ or $y' \lhd y$ then 
$\psi_\lambda(x'y') = \psi_\lambda(x')\psi_\lambda(y')$. In addition,  
$(\psi_\lambda(xy))(x_e) \not= (\psi_\lambda(x)\psi_\lambda(y))(x_e)$ can be supposed since the composition 
of continuous endomorphisms is continuous. Moreover, let $x,y$ be reduced. Otherwise there exist 
$h, \oh \in G$ such that $x\inv_\oh y$, $x\inv_h x x_\oh$ are reduced where 
$(\psi_\lambda((x\inv_h x x_\oh)(x\inv_\oh y)))(x_e) \not=
(\psi_\lambda(x\inv_h x x_\oh)\psi_\lambda(x\inv_\oh y))(x_e)$ by statement (7). Since 
$x\inv_h x x_\oh$ and $x\inv_\oh y$ have the same complexity as $x$ and $y$ respectively they are minimal 
in the same sense as $x$ and $y$. We discuss the three cases {\bf M1},{\bf M2}, and {\bf M3} as presented 
in the remarks about left representations of sums and products after Definition \ref{Definitionreduziert} 
using the corresponding notation and assumptions.\smallskip\\
{\bf M1:} Then $(\psi_\lambda(xy))(x_e) = \sum_{g \in \fC} (\psi_\lambda(xd_g)x_g)(x_e)$ where 
$\psi_\lambda(xd_g)x_g =\psi_\lambda(x)\psi_\lambda(d_g)x_g$ for all $g \in \fC$ since $x \unlhd x$ 
and $d_g \lhd y$. Now, the continuity of $\psi_\lambda(x)$ yields
\begin{eqnarray*}
(\psi_\lambda(x)\psi_\lambda(y))(x_e) &=& (\psi_\lambda(x))((\psi_\lambda(y))(x_e)) 
                                       =  (\psi_\lambda(x))(\sum_{g \in \fC} (\psi_\lambda(d_g)x_g)(x_e))\\ 
                                      &=& \sum_{g \in \fC} (\psi_\lambda(x))((\psi_\lambda(d_g)x_g)(x_e))
                                       = \sum_{g \in \fC} (\psi_\lambda(x)\psi_\lambda(d_g)x_g)(x_e)\\ 
                                      &=& \sum_{g \in \fC} (\psi_\lambda(xd_g)x_g)(x_e) = (\psi_\lambda(xy))(x_e).
\end{eqnarray*}
{\bf M2:} Then 
$(\psi_\lambda(xy))(x_e) = \sum_{g \in \fC} (\psi_\lambda(d_g \gap x_g \gap y \gap x\inv_g)\gap x_g)(x_e)$. 
Because of $d_g \lhd x$ and $y \unlhd y$ we obtain  
$\psi_\lambda(d_g \gap x_g \gap y \gap x\inv_g) \gap x_g =\psi_\lambda(d_g) \gap x_g \gap \psi_\lambda(y)$ 
for all $g \in \fC$. Moreover, $y \in D_C^-$ due to the specific situation of case $\bf M2$ and therefore 
$(\psi_\lambda(y))(x_e) \in F((C))$ by statement (6). Finally, statement (8) of 
Lemma \ref{aufbauvonpsiallgemeiner} provides
\begin{eqnarray*}
(\psi_\lambda(x)\psi_\lambda(y))(x_e) &=& (\psi_\lambda(x))((\psi_\lambda(y))(x_e)) 
                                       =  \sum_{g \in \fC} (\psi_\lambda(d_g) \gap x_g)((\psi_\lambda(y))(x_e))\\ 
                                      &=& \sum_{g \in \fC} (\psi_\lambda(d_g) \gap x_g \gap \psi_\lambda(y))(x_e)
                                       =  \sum_{g \in \fC} (\psi_\lambda(d_g \gap x_g \gap y \gap x\inv_g) \gap x_g)(x_e)\\ 
                                    &=& (\psi_\lambda(xy))(x_e).
\end{eqnarray*}

{\bf M3:} Then $(\psi_\lambda(xy))(x_e) = \sum_{g \in \fC} (\psi_\lambda(d_g) x_g)(x_e)$ 
with
\[d_g = \sum_{\og, \tg \in \fC \atop \og \tg C' = gC'} \od_\og \big( (x_\og \td_\tg x\inv_\og) x_\og x_\tg x\inv_g\big)
      = \sum_{\og, \tg \in \fC \atop \og \tg C' = gC'} \od_\og x_\og \td_\tg x_\tg x\inv_g \]
where $\od_\og \lhd x$ and $x_\og \td_\tg x_\tg x\inv_g \lhd y$ for all $\og,\tg \in \fC$. 
Because of statement (3) and the minimality property of $x, y$ regarding their complexities we conclude by 
arguments similar to those above that
\[\psi_\lambda(d_g) = \sum_{\og, \tg \in \fC \atop \og \tg C' = gC'} 
                      \psi_\lambda(\od_\og) x_\og \psi_\lambda(\td_\tg) x_\tg x\inv_g\]
for any $g \in \fC$ where all sums are finite. We now proceed as in case $\bf M2$ and apply for instance 
statement (8) of Lemma \ref{aufbauvonpsiallgemeiner} since $(\psi_\lambda(y))(x_e) \in F((C))$ 
by statement (6):
\begin{eqnarray*}
(\psi_\lambda(x)\psi_\lambda(y))(x_e) &=& \psi_\lambda(x)((\psi_\lambda(y))(x_e))\\
&=& \sum_{\og \in \fC}(\psi_\lambda(\od_\og) \gap x_\og)((\psi_\lambda(y))(x_e))\\
&=& \sum_{\og \in \fC}(\psi_\lambda(\od_\og) \gap x_\og)(\sum_{\tg \in \fC}(\psi_\lambda(\td_\tg) \gap x_\tg)(x_e))\\
&=& \sum_{\og \in \fC}\sum_{\tg \in \fC}(\psi_\lambda(\od_\og) x_\og)((\psi_\lambda(\td_\tg) x_\tg)(x_e))\\
&=& \sum_{\og \in \fC}\sum_{\tg \in \fC}(\psi_\lambda(\od_\og) \gap x_\og \gap \psi_\lambda(\td_\tg) \gap x_\tg)(x_e).
\end{eqnarray*}

For a fixed $\og \in \fC$ the sets 
$\supp (\psi_\lambda(\od_\og) \gap x_\og \gap \psi_\lambda(\td_\tg) \gap x_\tg)(x_e)$, $\tg \in \fC$ are 
pairwise disjoint since 
$(\psi_\lambda(\od_\og) \gap x_\og \gap \psi_\lambda(\td_\tg) \gap x_\tg)(x_e) \in x_\og x_\tg F((C'))$ 
by statement (6) and the induction hypothesis. Thus, Lemma \ref{summensummen} can be applied where    
$I = \{\og \in \fC \mid \od_\og \not= \Null\}$, $J = \{\tg \in \fC \mid \td_\tg \not= \Null\}$, and $K=\fC$ 
with $M_g = \{(\og, \tg) \in I \times J \mid \og \tg C' = gC'\}$ for all $g \in \fC$. We conclude
\begin{eqnarray*}
\sum_{\og \in \fC}\sum_{\tg \in \fC}(\psi_\lambda(\od_\og) \gap x_\og \gap \psi_\lambda(\td_\tg) \gap x_\tg)(x_e)&=&
\sum_{g \in \fC} \sum_{(\og,\tg) \in M_g} (\psi_\lambda(\od_\og) \gap x_\og \gap \psi_\lambda(\td_\tg) \gap x_\tg)(x_e)\\
&=& \sum_{g \in \fC} (\psi_\lambda(d_g) \gap x_g)(x_e) = (\psi_\lambda(xy))(x_e).
\end{eqnarray*}
In any case we deduce a contradiction which finally shows statement $\bf M$. So far, $\psi_\lambda(x)$ 
has been defined only if $x$ belongs to some $D_{\leq \lambda'}$ where $\lambda' < \lambda$. Now we 
treat the general situation and proceed by transfinite induction on the complexity $cp(x)$ using that 
$\psi_\lambda(x)$ has already been defined if $x \in F[G,\eta,\alpha]$ or $x \in D_{\leq \lambda'}$ for 
some $\lambda' < \lambda$. \smallskip\\
Case 1: $x$ is a proper atom. Then $x^{-1} \lhd x$ and $x^{-1} \in D_{\leq \lambda'}$ for some 
$\lambda' \in \Lambda$, $\lambda' < \lambda$ by \cite[Proposition 4.1]{DGH} and 
Proposition \ref{Dleqlambda}. According to the induction hypothesis it can be assumed 
that $\psi_\lambda(x^{-1})$ is a continuous and $v$-compatible automorphism of $F((G))$. We define 
$\psi_\lambda(x)=[\psi_\lambda(x^{-1})]^{-1}$ such that statement (5) follows by 
Proposition \ref{vvertraeglichistinjektiv} and Theorem \ref{finvstetig}. To prove statement (7) let 
$\tg \in G$. Then, $x_\tg x$ is a proper atom, $x_\tg x \in D_{\leq \lambda}$, and 
$\psi_\lambda(x_\tg x) = [\psi_\lambda((x_\tg x)^{-1})]^{-1} = [\psi_\lambda(x^{-1}x_\tg^{-1})]^{-1}
= [\psi_\lambda(x^{-1})x_\tg^{-1})]^{-1} = x_\tg [\psi_\lambda(x^{-1})]^{-1} = x_\tg \psi_\lambda(x)$. 
Clearly, $\psi_\lambda (x x_\tg) = \psi_\lambda (x) x_\tg$ can be derived similarly.\smallskip\\
We turn to the proof of statement (6). Let $x$ be reduced with pure left representation 
$x = \sum_{g \in \fC} d_g x_g$ in $D_C^-((C/C', \eta, \alpha))$ where $d_g \not=0$ for at least two 
$g \in \fC$. Then $x^{-1} = \sum_{g \in \fC} \od_g x_g$ in $D_C^-((C/C', \eta, \alpha))$ with 
$\od_g \not=0$ for at least two $g \in \fC$. This shows that $x\inv$ is reduced and that 
$x^{-1} = \sum_{g \in \fC}\od_g x_g$ is the pure left representation of $x^{-1}$. The sets 
$\supp (\psi_\lambda(d_g) x_g)(x_e)$ with $g \in \fC$ are pairwise disjoint 
since $\supp (\psi_\lambda(d_g) x_g)(x_e) \subseteq x_g F((C'))$ by induction hypothesis such that 
$\{(\psi_\lambda(d_g) x_g)(x_e) \mid g \in \fC\}$ is summable. We define 
$m:=\sum_{g \in \fC} (\psi_\lambda(d_g) x_g)(x_e) \in F((C))$. The arguments we shall apply 
now are essentially the same we already used before in subcase $\bf M3$ of the proof of statement $\bf M$. 
Hence, it will not be necessary to repeat all details again. We will use that $x\inv$ and each $d_g$, 
$\od_g$ belong to $D_{< \lambda'}$ with $\lambda' = cp(x)$ and that the restriction of $\psi_\lambda$ 
to $D_{< \lambda'}$ is an injective ring homomorphism. Moreover, Proposition \ref{Dleqlambda} provides 
$x\inv \in D_{\leq \lambda^\ast}$ for some $\lambda^\ast \in \Lambda$ with $\lambda^\ast < \lambda$ 
such that statements (5)-(8) of Lemma \ref{aufbauvonpsiallgemeiner} can be applied analogously to $x\inv$. 
We obtain
\begin{eqnarray*}
(\psi_\lambda(x^{-1}))(m) &=& \sum_{\og \in \fC}(\psi_\lambda(\od_\og) x_\og)(m) 
 =  \sum_{\og \in \fC}(\psi_\lambda(\od_\og) x_\og) (\sum_{g \in \fC}(\psi_\lambda(d_g) x_g)(x_e))\\
&=& \sum_{\og \in \fC}\sum_{g \in \fC}(\psi_\lambda(\od_\og x_\og d_g x_g))(x_e)
 =  \sum_{\tg \in \fC}\sum_{\og, g \in \fC \atop \og g C' = \tg C'} (\psi_\lambda(\od_\og x_\og d_g x_g))(x_e)\\
&=& \sum_{\tg \in \fC} (\psi_\lambda (\sum_{\og, g \in \fC \atop \og g C' = \tg C'} \od_\og x_\og d_g x_g x\inv_\tg)x_\tg)(x_e)
 =  \sum_{\tg \in \fC} (\psi_\lambda(\td_\tg) x_\tg)(x_e)
\end{eqnarray*}
where
\[\td_\tg = \sum_{\og, g \in \fC \atop \og g C' = \tg C'} \od_\og x_\og d_g x_g x\inv_\tg \mbox{\ for all $\tg \in \fC$\ }.\]
On the other hand $\sum_{\og \in \fC} \od_\og x_\og\gap\cdot\gap\sum_{g \in \fC} d_g x_g = \sum_{\tg \in \fC} \td_\tg x_\tg = 1x_e$ 
in $D^-_C((C/C', \eta, \alpha))$. Therefore, $(\psi_\lambda(x\inv))(m) = x_e$, that is
\[(\psi_\lambda(x))(x_e) = ([\psi_\lambda(x^{-1})]^{-1})(x_e) = m = 
\sum_{g \in \fC} (\psi_\lambda(d_g) x_g)(x_e) \in F((C)).\]
Case 2: $x$ is additively indecomposable. Let $x = x_1 \cdot \dots \cdot x_n$ be a complete multiplicative 
decomposition of $x$ which implies that any $x_i$ is a proper atom by \cite[Theorem 4.6]{DGH}. We define 
$\psi_\lambda (x)$ and prove statements (5),(6), and (7) by induction on $n$. The case $n=1$ 
has been treated above and for $n>1$ we introduce $\ox := x_2 \cdot \dots \cdot x_n$. Then, $\ox$ is 
additively indecomposable by \cite[Lemma 4.5]{DGH} with $\ox = x_2 \cdot \dots \cdot x_n$ as a complete 
multiplicative decomposition of $\ox$ by \cite[Theorem 4.6]{DGH}. Clearly, $x_1, \ox \lhd x$ yields 
$x_1, \ox \in D_{\leq \lambda}$ and according to Proposition \ref{Dleqlambda} there exists 
$\lambda' \in \Lambda$ with $\lambda' < \lambda$ such that $x_1, \ox \in D_{\leq \lambda'}$. Hence, 
$\psi_\lambda(x_1)$ and $\psi_\lambda(\ox)$ are given and we define 
$\psi_\lambda(x) = \psi_\lambda(x_1)\psi_\lambda(\ox)$. Later on we will show that this definition 
does not depend on the special decomposition $x = x_1 \cdot \dots \cdot x_n$. Clearly, statement (5) 
follows by Proposition \ref{speziellesstetig} and we shall prove statement (6). Thus, let $x$ be reduced.
There exist $h, \oh \in G$ such that $x_h^{-1}\ox$ and $x_\oh^{-1} x_1 x_h$ are reduced. 
\cite[Proposition 4.2]{DGH} ensures that $x = (x_1 x_h) \cdot (x_h^{-1}x_2) \cdot x_3 \cdot \dots \cdot x_n$ 
is a complete multiplicative decomposition of $x$ and 
$x_\oh^{-1} x = (x_\oh^{-1} x_1 x_h) \cdot (x_h^{-1} x_2) \cdot x_3 \cdot  \dots \cdot x_n$ is a complete 
multiplicative decomposition of  $x_\oh^{-1} x$ by Proposition \ref{neuezerlegungen}. We define 
$x' = x_\oh^{-1} x$, $x'_1 = x_\oh^{-1} x_1 x_h$, and $\ox' = x_h^{-1}\ox$. Now, let $a,b \in D$ arbitrary 
such that $ a\unlhd x'_1, b \unlhd \ox'$ whereby $a \lhd x'_1$ or $b \lhd \ox'$. Then, $a,b,ab \lhd x$ 
because of \cite[Theorem 4.6]{DGH} such that $\psi_\lambda(a), \psi_\lambda(b)$, and $\psi_\lambda(ab)$ are 
defined and moreover $\psi_\lambda(ab) = \psi_\lambda(a)\psi_\lambda(b)$. The latter statement is not 
obvious but it can be proved exactly the same way as $\psi_\lambda(xy) = \psi_\lambda(x) \psi_\lambda(y)$ 
in statement $\bf M$ and therefore an additional proof will be omitted.\smallskip\\
In order to show statement (6) we need the pure left representation of $x=x_1\ox$ which will be 
derived from suitable left representations of $x'_1$ or $\ox'$. We proceed as in the proof of 
statement (4) and again we have to distinguish three cases which correspond to the three cases {\bf M1}, 
{\bf M2}, and {\bf M3} as discussed in the remarks about left representations of sums and
products after Definition 6.3.\smallskip\\
{\bf M1:} Then $\ox' = \sum_{g \in \fC} d_g x_g$ and $x'_1 \in D^-_C$ such that 
$x_\oh^{-1} x = x' = \sum_{g \in \fC} (x'_1 d_g) x_g$ is a left representation of $x$ and $x'$ where 
$\oh \in C$ since $x$ is reduced. For any $g \in \fC$ there exists some $g' \in \fC$ satisfying 
$\oh g C' = g'C'$, that is $x_\oh x_g x\inv_{g'}\in F^\times X_{C'}$. This yields 
$x = \sum_{g \in \fC} ((x_\oh \gap x'_1 \gap d_g x\inv_\oh)(x_\oh x_g x\inv_{g'})) \gap x_{g'}$ 
as the pure left representation of $x$. Applying the induction hypothesis to $\ox'$ and 
Lemma \ref{aufbauvonpsiallgemeiner} we obtain by means of $d_g \lhd \ox'$ for all $g \in \fC$:
\begin{eqnarray*}
(\psi_\lambda(x))(x_e) &=& (\psi_\lambda(x_1))((\psi_\lambda(\ox))(x_e))\\ 
                     &=& (\psi_\lambda(x_\oh x'_1 x_h^{-1}))((\psi_\lambda(x_h\ox'))(x_e))\\
                     &=& x_\oh(\psi_\lambda(x'_1))((\psi_\lambda(\ox'))(x_e))\\ 
                     &=& x_\oh(\psi_\lambda(x'_1)(\sum_{g \in \fC} (\psi_\lambda(d_g) x_g)(x_e))\\
                     &=& x_\oh\sum_{g \in \fC} (\psi_\lambda(x'_1)\psi_\lambda(d_g) x_g)(x_e)\\
                     &=& x_\oh\sum_{g \in \fC} (\psi_\lambda(x'_1 d_g)x_g)(x_e)\\
                     &=& \sum_{g \in \fC} (\psi_\lambda ((x_\oh \gap x'_1 \gap d_g x\inv_\oh)(x_\oh x_g x\inv_{g'})) x_{g'})(x_e) \in F((C)).
\end{eqnarray*}
Thus, statement (6) is shown.\smallskip\\
{\bf M2:} Then $x_1' = \sum_{g \in \fC} d_g x_g$ and $\ox' \in D^-_C$ such that 
$x_\oh^{-1} x = x' = \sum_{g \in \fC} (d_g \gap x_g \gap \ox' \gap x_g^{-1}) \gap x_g$ is a left 
representation of $x$ and $x'$. Again, $\oh \in C$ since $x$ is reduced. The same arguments we just 
applied now show that 
$x = \sum_{g \in \fC} ((x_\oh \gap d_g \gap x_g \gap \ox' \gap x_g^{-1} x\inv_\oh)( x_\oh x_g x\inv_{g'})) \gap x_{g'}$ 
is the pure left representation of $x$ where $x_\oh x_g x\inv_{g'}\in F^\times X_{C'} $ for all $g \in \fC$. 
Arguments we repeatedly used before finally provide
\begin{eqnarray*}
(\psi_\lambda(x))(x_e) &=& (\psi_\lambda(x_1))((\psi_\lambda(\ox))(x_e)) 
  = x_\oh (\psi_\lambda(x'_1))((\psi_\lambda(\ox'))(x_e))\\
 &=& x_\oh \sum_{g \in \fC} (\psi_\lambda(d_g) x_g)((\psi_\lambda(\ox'))(x_e))
  =  x_\oh \sum_{g \in \fC} (\psi_\lambda(d_g x_g \ox' x_g^{-1}) x_g)(x_e)\\
 &=& \sum_{g \in \fC} (\psi_\lambda ((x_\oh \gap d_g \gap x_g \gap \ox' \gap x_g^{-1} x\inv_\oh) (x_\oh x_g x\inv_{g'})) \gap x_{g'})(x_e) 
     \in F((C))
\end{eqnarray*}
and therefore statement (6).\smallskip\\
{\bf M3:} Then, $x_1' = \sum_{\og \in \fC} \od_\og x_\og$ and $\ox' = \sum_{\tg \in \fC} \td_\tg x_\tg$ 
are the pure left representations of $x'_1$ and $\ox'$ respectively such that  
$x_\oh^{-1} x = x' = \sum_{g \in \fC} d_g x_g$ is a left representation of $x$ and $x'$ where $\oh \in C$ 
since $x$ is reduced and
\[d_g = \sum_{\og, \tg \in \fC \atop \og \tg C' = g C'} \od_\og x_\og \td_\tg x_\tg x\inv_g \mbox{\ for all $g \in \fC$\ }.\]

Therefore, $x = \sum_{g \in \fC} ((x_\oh d_g x\inv_\oh) (x_\oh x_g x\inv_{g'})) x_{g'}$ is the pure left  
representation of $x$ where $g' \in \fC$ and $hgC' = g'C'$ for all $g \in \fC$. We obtain
\begin{eqnarray*}
(\psi_\lambda(x))(x_e) &=& (\psi_\lambda(x_1))((\psi_\lambda(\ox))(x_e))\\
 &=& x_\oh(\psi_\lambda(x'_1))((\psi_\lambda(\ox'))(x_e))\\ 
 &=& x_\oh \sum_{\og \in \fC} (\psi_\lambda(\od_\og) x_\og)((\psi_\lambda(\ox'))(x_e))\\
 &=& x_\oh \sum_{\og \in \fC} \sum_{\tg \in \fC}(\psi_\lambda(\od_\og)x_\og)((\psi_\lambda(\td_\tg)x_\tg)(x_e))\\
 &=& x_\oh \sum_{\og \in \fC} \sum_{\tg \in \fC}(\psi_\lambda(\od_\og x_\og \td_\tg x_\tg))(x_e)\\
 &=& x_\oh \sum_{g \in \fC} \sum_{\og, \tg \in \fC \atop \og \tg C' = g C'} (\psi_\lambda(\od_\og x_\og \td_\tg x_\tg x\inv_g) x_g)(x_e)\\
 &=& x_\oh \sum_{g \in \fC} (\psi_\lambda (d_g) x_g)(x_e)\\
 &=& \sum_{g \in \fC} (\psi_\lambda(x_\oh d_g x\inv_\oh x_\oh x_g x\inv_{g'}) x_{g'})(x_e) \in F((C))
\end{eqnarray*}
which shows statement (6).\smallskip\\
We still need to prove that $\psi_\lambda(x)$ is well-defined. Thus, let 
$x = x_1 \cdot x_2 \cdot \dots \cdot x_n$, $x  = y_1 \cdot y_2 \cdot \dots \cdot y_m$ be two complete 
multiplicative decompositions of $x$. We define $\ox = x_2 \cdot \dots \cdot x_n$ and  
$\oy = y_2 \cdot \dots \cdot y_{m}$. Since both endomorphisms are continuous it is enough to show 
$(\psi_\lambda(x_1)\psi_\lambda(\ox))(x_g)=(\psi_\lambda(y_1)\psi_\lambda(\oy))(x_g)$, that is, 
$(\psi_\lambda(x_1)\psi_\lambda(\ox\gap x_g))(x_e) = (\psi_\lambda(y_1)\psi_\lambda(\oy\gap x_g))(x_e)$ 
for an arbitrary $g \in G$. According to (the remark after) Proposition \ref{neuezerlegungen} the element 
$x x_g$ is additively indecomposable with two complete additive decompositions 
\[x x_g = x_1 \cdot \dots \cdot x_{n-1} \cdot (x_n x_g) = y_1 \cdot \dots \cdot y_{m-1} \cdot (y_m x_g).\]
Without loss of generality, we may therefore suppose $g=e$. Since $x\inv _h x$ is reduced for some 
$h \in G$ we may also assume for the same reason that $x$ itself is reduced. But then 
$(\psi_\lambda(x))(x_e)$ depends only on the pure left-representation of $x$ and not on the complete 
multiplicative decomposition as we have just proved. This shows, that $\psi_\lambda(x)$ is well-defined 
if $x$ is additively indecomposable and also statement (6) is shown. We finally turn to statement (7) 
and consider $\tg \in G$. Clearly, $x_\tg x = (x_\tg x_1) \cdot x_2 \cdot \dots \cdot x_n$ is a complete 
multiplicative decomposition of the additively indecomposable element $x_\tg x$ by Proposition 
\ref{neuezerlegungen}. Hence,
\begin{eqnarray*}
\psi_\lambda(x_\tg x)&=&\psi_\lambda(x_\tg x_1)\cdot\psi_\lambda(x_2 \cdot\dots\cdot x_n)\\
                   &=&x_\tg \psi_\lambda(x_1)\cdot\psi_\lambda(x_2 \cdot\dots\cdot x_n) = x_\tg\psi_\lambda(x)
\end{eqnarray*}
and the claim for $x x_\tg$ follows similarly.\smallskip\\
Case 3: Now, let $x$ be an arbitrary element having a complete additive decomposition $x = x_1 + \dots + x_n$, 
$n \in \N$. Then any $x_i$, $i = 1,\dots,n$  is additively indecomposable by \cite[Theorem 3.6]{DGH}. 
We define $\psi_\lambda(x)$ and prove statemants (6) and (7) by induction on $n$. The case $n=1$ has 
been treated above. Thus, let $n>1$ and $\ox = x_2 + \dots + x_n$ which presents in addition 
a complete additive decomposition of $\ox$ due to \cite[Theorem 3.6]{DGH}. Clearly, $x_1, \ox \lhd x$ 
yields $x_1, \ox \in D_{\leq \lambda}$ by Proposition \ref{Dleqlambda} such that $\psi_\lambda(x_1)$, 
$\psi_\lambda(\ox)$ are already given by induction hypothesis. We define 
$\psi_\lambda(x) = \psi_\lambda(x_1) + \psi_\lambda(\ox)$ and will prove later on that this definition 
does not depend on the complete additive decomposition of $x$ as introduced above. In order to show 
statement (6) let $x$ be reduced and let $cp(x_1) > 1$ or $cp(\ox) > 1$. The pure left representation 
of $x = x_1 + \ox$ results from left representations of $x_1$ or $\ox$ where two different 
situations must be considered which correspond to the cases $\bf A1$ and $\bf A2$ .\smallskip\\
{\bf A1:} Without loss of generality (if necessary the roles of $x_1$ and $\ox$ have to be exchanged) there 
exist $h\in G$, $C \in \Co^\ast_\leq$, and $g \in \fC, g\not=e$ such that 
$x\inv_h x_1, x\inv_h \ox x\inv_g \in D^-_C$ and
\[x_h^{-1} x = x_h^{-1}(x_1 + \ox) = x_h^{-1} x_1 + (x_h^{-1} \ox x_g^{-1}) x_g\]
is a left representation of $x$. Since $x$ is reduced $h \in C$ follows and  
$x_1 x_h\inv, \ox x_g^{-1} x_h\inv \in D^-_C$. There exist $\oh, \og \in \fC$ satisfying $hC' = \oh C'$ and 
$hgC' = \og C'$ where $\oh \not= \og$ because of $g \not\in C$. Hence 
$x_1 x\inv_\oh = (x_1 x\inv_h) (x_h x\inv_\oh), \ox x\inv_\og = (\ox x\inv_g x\inv_h) (x_h x_g x\inv_\og) \in D^-_C$ 
such that 
\[x = (x_1 x\inv_\oh)x_\oh + (\ox x\inv_\og)x_\og\]
is the pure left representation of $x$ and statement (6) follows from 
\begin{eqnarray*}
(\psi_\lambda(x))(x_e) &=& (\psi_\lambda(x_1))(x_e) + (\psi_\lambda(\ox))(x_e)\\
                     &=& (\psi_\lambda(x_1 x\inv_\oh)x_\oh)(x_e) + (\psi_\lambda(\ox x\inv_\og)x_\og)(x_e).
\end{eqnarray*}
{\bf A2:} Then, there exist left representations $x_h^{-1}x_1 = \sum_{g \in \fC} d_g x_g$ and 
$x_h^{-1}\ox = \sum_{g \in \fC} d_g' x_g$ for $x_1$ and $\ox$ respectively. Because of 
$d_g + d_g' \lhd x_1 + \ox = x$ the subcase $\bf A2.2$ can not occur such that 
$x_h^{-1} x = x_h^{-1} x_1 + x_h^{-1} \ox = \sum_{g \in \fC} (d_g + d_g') x_g$ is a left representation 
of $x$. Since $x$ is reduced we obtain
\[x = \sum_{g \in \fC} ((x_h(d_g + d_g')x\inv_h)(x_h x_g x\inv_{g'}))x_{g'}\]
as the pure left representation of $x$ where $g' \in \fC$ and $hgC' = g'C'$ for all $g \in \fC$. Now, 
statement (6) follows from
\begin{eqnarray*}
(\psi_\lambda(x))(x_e) &=& (\psi_\lambda(x_1))(x_e) + (\psi_\lambda(\ox))(x_e)\\
                       &=& x_h \sum_{g \in \fC} [(\psi_\lambda(d_g) x_g)(x_e) + (\psi_\lambda(d'_g) x_g)(x_e)]\\
                       &=& x_h \sum_{g \in \fC} ((\psi_\lambda(d_g) + \psi_\lambda(d'_g)) x_g)(x_e)\\
                       &=& x_h \sum_{g \in \fC} (\psi_\lambda(d_g + d'_g) x_g)(x_e)\\
                       &=& \sum_{g \in \fC} (\psi_\lambda ((x_h(d_g + d_g')x\inv_h)(x_h x_g x\inv_{g'}))x_{g'})(x_e).
\end{eqnarray*}
To show that $\psi_\lambda(x)$ is well-defined let 
$x = x_1 + \dots + x_n$, $x = y_1 + \dots + y_m$ be two complete additive decompositions of $x$ and 
$\ox = x_2 + \dots + x_n$,  $\oy = y_2 + \dots + y_{m}$. We assume 
$\psi_\lambda(x_1) + \psi_\lambda(\ox) \not= \psi_\lambda(y_1) + \psi_\lambda(\oy)$, that is, 
$(\psi_\lambda(x_1) + \psi_\lambda(\ox))(x_g) \not= (\psi_\lambda(y_1) + \psi_\lambda(\oy))(x_g)$ and 
therefore $(\psi_\lambda(x_1 \gap x_g) + \psi_\lambda(\ox \gap x_g))(x_e) \not= 
(\psi_\lambda(y_1 \gap x_g) + \psi_\lambda(\oy \gap x_g))(x_e)$ for some $g \in G$ since both endomorphisms 
are continuous. Furthermore, $x x_g = x_1 x_g + \dots + x_n x_g$, $x x_g = y_1 x_g + \dots + y_m x_g$ are 
two complete additive decompositions of $x x_g$ by Proposition \ref{neuezerlegungen} which means that we 
may restrict to the case $g=e$. Since $x\inv _h x$ is reduced for some $h \in G$ we may also assume for 
the same reason that $x$ itself is reduced. But then $(\psi_\lambda(x))(x_e)$ depends only on the pure 
left-representation of $x$ and not on the complete additive decomposition as we have just proved. This 
shows, that $\psi_\lambda(x)$ is well-defined if $x$ is additively indecomposable and also statement (6) 
is shown. In order to verify statement (7) we consider an arbitrary $\tg \in G$. Again, 
$x_\tg x = x_\tg x_1 + x_\tg x_2 + \dots + x_\tg x_n$ is a complete additive decomposition of $x_\tg x$ 
by Proposition \ref{neuezerlegungen} such that
\begin{eqnarray*}
\psi_\lambda(x_\tg x)&=&\psi_\lambda(x_\tg x_1) + \psi_\lambda(x_\tg x_2 + \dots + x_\tg x_n)\\
                     &=&x_\tg \psi_\lambda(x_1) + x_\tg \psi_\lambda(x_2 + \dots + x_n) = x_\tg\psi_\lambda(x).
\end{eqnarray*}
Similar arguments provide the claim for $x x_\tg$.\smallskip\\
It remains to prove statement (5). Since $\psi_\lambda(x)$ is continuous by induction hypothesis and 
Proposition \ref{speziellesstetig} we apply Theorem \ref{lokalglobal} where $x \not=0$ and show that 
$\psi_\lambda(x)$ is surjective on $G$ and $v$-compatible on $G$.\smallskip\\
$\psi_\lambda(x)$ is $v$-compatible on $G$: We assume that $\psi_\lambda(y)$ is $v$-compatible on $G$ 
for each $y \in D$, $y \not= 0$ where $y \lhd x$ and consider $\og, \tg \in G$ with $\og < \tg$. Then, 
$v((\psi_\lambda(x))(x_\og)) < v((\psi_\lambda(x))(x_\tg))$ has to be shown, that is, 
$v((\psi_\lambda(x x_\og))(x_e)) < v((\psi_\lambda(x x_\og))(x_\og^{-1} x_\tg)) = 
v((\psi_\lambda(x x_\og))(x_{\og^{-1} \tg}))$ where $e < \og^{-1} \tg$. Therefore, without loss of 
generality let us assume $\og = e$. Moreover, we may restrict to the case where $cp(x) > 1$ and $x$ is 
reduced. Let $x = \sum_{g \in \fC} d_g x_g \in D^-_C((C/C', \eta, \alpha))$ be the pure left representation 
of $x$ and where $g_0$ is the minimal element in the support of $x$. Two cases have to be differentiated 
here:\smallskip\\
Case 1: $\tg \not\in C$. Then $C \subseteq C_\tg^-$ and $x \in D^+_C \subseteq D_\tg^-$. 
This means $x_\tg^{-1} x x_\tg \in D_\tg^-$, that is, $(\psi_\lambda(x_\tg^{-1} x x_\tg))(x_e) \in F((C_\tg^-))$ 
and $(\psi_\lambda(x))(x_\tg) = x_\tg(\psi_\lambda(x_\tg^{-1} x x_\tg))(x_e) \in x_\tg F((C_\tg^-))$ 
follows. Because of $(\psi_\lambda(x))(x_e) \in F((C)) \subseteq F((C_\tg^-))$ and $\tg > e$ we are 
done.\smallskip\\
Case 2: $\tg \in C$. For $g \in \fC$ there exists $g' \in \fC$ satisfying $g \tg C' = g'C'$. Since $C/C'$ 
is abelian we conclude $g' < \og'$ for all $g, \og \in \fC$ with $g < \og$. This shows that 
\[x x_\tg = \sum_{g \in \fC} (d_g (x_g x_\tg x\inv_{g'}))x_{g'}\]
is a left representation of $x x_\tg$ where $g_0'$ is the minimal element in the support of $x x_\tg$. 
For all $g \in \fC$ we obtain 
\[(\psi_\lambda(d_g (x_g x_\tg x\inv_{g'}))x_{g'})(x_e) \in x_{g'}F((C')),\]
that is
\[v((\psi_\lambda(x))(x_\tg)) = v((\psi_\lambda(x x_\tg))(x_e))
                              = v((\psi_\lambda(d_{g_0} x_{g_0} x_\tg))(x_e))
                             = v((\psi_\lambda(d_{g_0} x_{g_0}))(x_\tg)).\]
Clearly, $v((\psi_\lambda(x))(x_e))= v((\psi_\lambda(d_{g_0} x_{g_0}))(x_e))$ follows in a similar 
fashion. Because of $d_{g_0}\lhd x$ and $\tg > e$ we obtain 
$v((\psi_\lambda(d_{g_0} x_{g_0}))(x_e)) < v((\psi_\lambda(d_{g_0} x_{g_0}))(x_\tg))$ as claimed.
\smallskip\\
$\psi_\lambda(x)$ is surjective on $G$: We assume that $\psi_\lambda(y)$ is surjective on $G$ for all 
$y \in D$, $y \not= 0$ where $y \lhd x$ and choose an arbitrary $\tg \in G$. It has to be shown that 
$v((\psi_\lambda(x))(x_\og))= \tg$ holds true for some $\og \in G$ which is equivalent to 
$v((\psi_\lambda(x_\tg^{-1} x))(x_\og))= e$. Thus, let $\tg = e$ without loss of generality and let 
$x x_h^{-1} = \sum_{g \in \fC} d_g x_g$ be a right representation of $x$. If there is some $\og \in G$ 
satisfying $v((\psi_\lambda(x x_h^{-1}))(x_\og))= e$ then $v((\psi_\lambda(x))(x_{h^{-1}\og}))= e$ 
and vice versa. Therefore, we may choose $h = e$, that is $x = \sum_{g \in \fC} d_g x_g$. Because of 
$d_{g_0} \lhd x$ there exists some $\og \in G$ satisfying 
\[v((\psi_\lambda(d_{g_0} x_{g_0} x_\og))(x_e)) = v((\psi_\lambda(d_{g_0} x_{g_0}))(x_\og)) = e\]
where $g_0$ is minimal in the support of $x$. If $\og \not\in C$, which means $C \subseteq C_\og^-$, 
then $D_C^+ \subseteq D_\og^-$ and $x_\og^{-1} d_{g_0} x_{g_0} x_\og \in D_\og^-$ since 
$d_{g_0} x_{g_0} \in D^-_C \subseteq D_\og^-$. Hence, 
$x_\og (\psi_\lambda(x_\og^{-1} d_{g_0} x_{g_0} x_\og))(x_e)$ is in $x_\og F((C_\og^-))$. But this would 
be contrary to $v((\psi_\lambda(d_{g_0} x_{g_0} x_\og))(x_e)) = e$. Thus,  
$\og \in C$ follows and again for any $g \in \fC$ let $g'$ be in $\fC$ satisfying $g \og C' = g' C'$. 
Similar to the arguments we used above we conclude 
$(\psi_\lambda(x))(x_\og) = \sum_{g \in \fC}(\psi_\lambda(d_g\gap (x_g x_\og x\inv_{g'})) x_{g'})(x_e)$ 
which provides 
\[v((\psi_\lambda(x))(x_\og)) = v((\psi_\lambda(d_{g_0}\gap x_{g_0}) x_\og)(x_e))
                              = v((\psi_\lambda(d_{g_0}\gap x_{g_0}))( x_\og)) = e.\]
\qed\\

We finally turn to the proof of Theorem \ref{Fortsetzungvonphi}. For any $x \in D$ there exists some 
$\lambda \in \Lambda$ satisfying $x \in D_{\leq \lambda}$ and Lemma \ref{aufbauvonpsi} now shows that
$\psi(x): = \psi_\lambda(x)$ where $\lambda \in \Lambda$ and $x \in D_{\leq\lambda}$ yields a well-defined 
mapping $\psi: D \longrightarrow \R$. In order to show that $\psi$ is additive let us consider arbitrary 
elements $x,y \in D$. Then $x,y \in D_{\leq\lambda}$ for some $\lambda \in \Lambda$. We have to verify 
$\psi_\lambda(x+y)=\psi_\lambda(x)+\psi_\lambda(y)$ which can be done as in the proof of statement 
$\bf A$. Using the notation introduced there let us assume $\psi_\lambda(x+y)\not=\psi_\lambda(x)+\psi_\lambda(y)$ 
and that $(m^+,m^-)$ is minimal with respect to the lexicographical order among all $x,y \in D_{\leq\lambda}$. 
Now, the arguments we applied in the proof of statement $\bf A$ can be adopted straightforwardly. Similarly, 
the multiplicativity of $\psi$ can be derived as statement $\bf M$. Thus, $\psi$ is a ring homomorphism which 
in addition is injective because of statement (5) of Lemma \ref{aufbauvonpsi}. Therefore, $\psi(D)$ is a 
subdivision ring of $\R$ containing $F[G,\eta,\alpha]$ as a subring. Since $\R$ is a division ring of fractions 
of $F[G,\eta,\alpha]$ we conclude that $\psi$ is surjective.

\section{Hughes' Theorems and an Example}\label{Hughes Theorems}

\vspace*{0.3cm}

We now discuss I. Hughes' two main theorems from \cite{H1,H2} in the light of our results.

\begin{theorem}\label{Hughes1}
Let $R = F[G,\eta,\alpha]$ be a crossed product of a locally indicable group $G$ over a skew field 
$F$. If $R$ has two Hughes-free division rings of fractions $D$ and $D'$, then they 
are isomorphic by an isomorphism which extends the identity map on $R$.
\end{theorem}

Proofs are given in \cite{H1,DHS,J}. According to the arguments from \cite{J} we obtain this 
theorem from Corollary \ref{alleisomorph} and the fact that both division rings are free with respect 
to a Conradian left-order (with maximal rank) of $G$. Following \cite{H2} we call a locally 
indicable group $G$ freely embeddable if  $F[G,\eta,\alpha]$ has a Hughes-free division ring of 
fractions for each skew field $F$. The next theorem is the main result of \cite{H2}.

\begin{theorem}\label{Hughes2}
Let $G$ be a locally indicable group and let $N$ be a normal subgroup of $G$. If $G/N$ and $N$ are 
both freely embeddable then $G$ is also freely embeddable.
\end{theorem}

{\bf Proof.} Let $F[G,\eta,\alpha]$ be an arbitrary crossed product of $G$ over a skew field $F$. 
According to Proposition \ref{RGverschraenktueberRN} we consider $F[G,\eta,\alpha]$ as a crossed 
product of $G/N$ over the ring $F[N,\eta,\alpha]$. Therefore, we shall write 
$(F[N,\eta,\alpha])[G/N,\oeta,\oalpha]$ instead of $F[G,\eta,\alpha]$ where $\oalpha_g$ is the 
restriction of $\alpha_g$ to $F[N,\eta,\alpha]$ for all $g \in G$ which is an automorphism 
of $F[N,\eta,\alpha]$. By assumption, $F[N,\eta,\alpha]$ possesses a Hughes-free division ring of 
fractions $D_N$. Because of Theorem \ref{Hughes1} (or Corollary \ref{alleisomorph}) any $\oalpha_g$, 
$g \in G$ can be extended to an automorphism of $D_N$. This extension is unique since $D_N$ is a 
division ring of fractions of $F[N,\eta,\alpha]$ and therefore it will also be denoted by $\oalpha_g$. 
As explained in Section \ref{Verschraenkte Produkte} before Proposition \ref{verschraenktesproduktueberd} 
the crossed product $D_N[G/N,\oeta,\oalpha]$ is well-defined and contains $F[G,\eta,\alpha]$ as a subring 
in the obvious meaning. By assumption, $D_N[G/N,\oeta,\oalpha]$ possesses a Hughes-free division ring of 
fractions $D$. Clearly, $D$ is a division ring of fractions of $F[G,\eta,\alpha]$ and it remains to prove 
that it is Hughes-free. We show that $D$ is free with respect to a Conradian left-order $\leq$ of $G$. 
Let $\leq_1$ and $\leq_2$ be Conradian left-orders of $N$ and $G/N$ respectively. Then, 
$P = P_{\leq_1} \cup \{g \in G \mid eN <_2 gN\}$ is a positive cone of $G$ of a suitable left-order 
$\leq$ of $G$. A subgroup $U$ of $G$ is convex in $G$ with respect to $\leq$ if and only if 
$U \subseteq N$ and $U$ is convex in $N$ with respect to $\leq_1$ or $N \subseteq U$ and $U/N$ is a 
convex subgroup of $G/N$ with respect to $\leq_2$. This shows that $\leq$ is Conradian. We consider 
an arbitrary convex jump $(C',C)$ of $G$ with respect to $\leq$ and arbitrarily given 
$h_1, \dots ,h_n \in C$ such that $h_1C', \dots , h_nC'$ are pairwise different. Then 
$C' \subset C \subseteq N$ or $N \subseteq C' \subset C$. In the first case we obtain that $(C',C)$ 
is a convex jump of $N$ and that the rational closure of $F[C',\eta,\alpha]$ in $D$ coincides with the 
rational closure $D_C^-$ of $F[C',\eta,\alpha]$ in $D_N$. Since $D_N$ is a free division ring of fractions 
of $F[N,\eta,\alpha]$ with respect to $\leq_1$ we conclude that $x_{h_1}, \dots, x_{h_n}$ are linearly 
independent over $D_C^-$. In the second case we obtain that $D_C^-$ is the rational closure of 
$F[C',\eta,\alpha]$ in $D$. Then $(C'/N, C/N)$ is a convex jump of $G/N$ with respect to $\leq_2$ by 
definition of $\leq$ and $(h_1N)C'/N, \dots , (h_nN) C'/N$ are pairwise different. Moreover, $D_C^-$ 
is the rational closure of $D_N[C'/N, \oeta,\oalpha]$ in $D$. Applying the notation of Proposition 
\ref{RGverschraenktueberRN} we choose $h'_1, \dots h'_n \in \frak G$ satisfying 
$h'_1 N = h_1 N, \dots, h'_n N = h_n N$. Then, $x_{h'_1}, \dots , x_{h'_n}$ are linearly independent 
over $D_C^-$ since $D$ is a free division ring of fractions of $D_N[G/N,\oeta,\oalpha]$ with respect 
to $\leq_2$ and the same holds true for $x_{h_1}, \dots , x_{h_n}$.\qed\\

We finally apply our results to a locally indicable group $G$ with a left-order which is not Conradian 
and refer the reader to \cite{D1,D2,GS} for further details. Let $G$ be the knot group of the trefoil 
which is isomorphic to the braid group $B_3$ on $3$ strands. It can be written as a semi-direct product 
$B_3 \cong \mbox{$F_2 \sdpr_\varphi \langle x \rangle $}$ where $F_2$ is the free group of rank $2$ 
and $\langle x \rangle$ an infinite cyclic group. As a free group $F_2$ is endowed with a two-sided 
order and therefore the group ring $F[F_2]$ of $F_2$ over an arbitrary skew field $F$ possesses a 
Hughes-free division ring of fractions $D_{F_2}$ which is given, for instance, as the rational closure 
of $F[F_2]$ in the corresponding Mal'cev-Neumann division ring of all formal power series in $F_2$ over 
$F$. As explained in the proof of Theorem \ref{Hughes2} the group ring $F[B_3]$ occurs as a subring of 
the skew polynomial ring $D_{F_2}[x, \varphi]$ which is a left and right Ore-domain. The corresponding 
division ring of fractions $D_{F_2}(x, \varphi)$ is a Hughes-free division ring of fractions of $F[B_3]$. 
On the other hand $B_3$ can be written as $B_3 = \langle u,w \mid wu^2w = u \rangle$ where $u$ and $w$ 
generate a semigroup $P \subseteq B_3$ which is a positive cone of a left-order $\leq'$ of $B_3$ that 
is not Conradian. The corresponding Dubrovin-quotient ring $\R$ of $F[B_3]$, that is, the rational 
closure of $F[B_3]$ in the endomorphism ring of the right $F$-vector space $F((B_3))$ of all formal 
power series in $B_3$ over $F$ with respect to $\leq'$ is a division ring. Even though $\leq'$ is not 
Conradian the division ring $\R$ is Hughes-free by Corollary \ref{dubrovinisthughes}, that is, $\R$ and 
$D_{F_2}(x, \varphi)$ are isomorphic. Therefore, $D_{F_2}(x, \varphi)$ has a rank 1 valuation ring 
possessing a prime ideal which is not completely prime.


%
\end{document}